\documentclass[12pt]{amsart}

\usepackage{amssymb}
\usepackage{amsmath}
\usepackage{amscd}
\usepackage{float}
\usepackage{graphicx}
\usepackage{amsfonts}
\usepackage{pb-diagram}
\usepackage{eufrak}

\newtheorem{thm}{Theorem}

\newtheorem{lem}{Lemma}

\newtheorem{defn}{Definition}
\newtheorem{rem}{Remark}
\newtheorem{conj}{Conjecture}

\begin{document}

\title[Braid equivalences and the $L$--moves]{Braid equivalences and the $L$--moves}

\author{Sofia Lambropoulou}
\address{ Department of Mathematics,
National Technical University of Athens,
Zografou Campus, GR-157 80 Athens, Greece.}
\email{sofia@math.ntua.gr}
\urladdr{http://www.math.ntua.gr/$\tilde{~}$sofia}

\keywords{Braiding, Markov theorem, $L$--moves, conjugation, commuting, stabilization, links in knot complements, links in $3$--manifolds, mixed links, mixed braids, braid equivalence, twisted loop conjugation, twisted band move, links in handlebodies, virtual braids, welded braids, singular braids.}

\subjclass[2000]{57M25, 57M27, 57N10}

\date{}
\maketitle

\begin{abstract}

In this survey paper we present the $L$--moves between braids and how they can adapt and serve for establishing and proving braid equivalence theorems for various diagrammatic settings, such as for classical knots, for knots in knot complements, in c.c.o. $3$--manifolds and in handlebodies, as well as for virtual knots, for flat virtuals, for welded knots and for singular knots. The $L$--moves are local and they provide a uniform ground for formulating and proving braid equivalence theorems for any diagrammatic setting where the notion of braid and diagrammatic isotopy is defined, the statements being first geometric and then algebraic.

\end{abstract}

\section*{Introduction}

 The central problem in classical knot theory is the complete classification of knots and links up to ambient isotopy in three-space. Reidemeister \cite{Rd1} translated isotopy in space into a diagrammatic equivalence relation.
Further, the braid group provides a fundamental algebraic structure associated with knots. The Alexander theorem tells us that every knot or link can be isotoped to braid form. The capstone of this relationship is the  Markov theorem, giving necessary and sufficient conditions for two braids to close to isotopic links. Since the discovery of the Jones polynomial \cite{Jo} --the first invariant for knots and links that was constructed using braids-- the Markov
theorem received renewed attention.

\smallbreak
Among other proofs of the Markov theorem, in
\cite{LR1} we introduced a new type of braid move, the $L$--move, and we proved an {\it one-move}
Markov theorem. The $L$--move is a very simple uniform geometric move that can be applied anywhere in a braid to produce a braid with the same closure. More precisely, an {\it $L$--move} consists in cutting a strand of the braid and taking the top end to the bottom of the braid and the bottom end to the top of the braid, both entirely over or both entirely under the braid, creating a new pair of corresponding braid strands.  By small braid isotopies, an $L$--move is equivalent to adding an in-box crossing in the braid, aquiring thus an algebraic expression.  See Fig.~5 for illustrations.
 Consequently, in \cite{LR1,LR2,HL} the $L$--move approach was used for proving braid equivalence theorems (analogues of the Markov theorem) for classical knots and links in knot complements, in closed, connected, oriented (c.c.o.) $3$--manifolds and in handlebodies. Further, using the $L$--move methods, we proved in \cite{KL2} braid equivalence theorems for the virtual braid group, for the welded braid group and for some analogues of these structures. Finally, in \cite{La3}  an $L$--move singular braid equivalence is given, using the algebraic singular braid equivalence of  \cite{Ge}.

\smallbreak

Given a diagrammatic knot theory there are deep interrelations between the diagrammatic knot isotopy in this theory, the braid structures and the corresponding braid equivalence. On the other hand, braid equivalence theorems are important for understanding the structure and classification of knots and links in various settings. Also, for constructing invariants of knots and links using algebraic means.

\smallbreak

In this survey paper we present the $L$--move methods of the author and collaborators and how they can serve for establishing and proving braid equivalence theorems for the diagrammatic settings mentioned above. In Section~1 we discuss braiding and braid equivalence for knots in $S^3$. We describe a generic braiding method that can adapt to any diagrammatic setting, we define the $L$--moves, and we give the one-move Markov theorem  (Theorem~\ref{lthm}). In Section~2 we present the $L$--move methods for deriving braid equivalences, geometric and algebraic, for knots in knot complements, in c.c.o. $3$--manifolds and in handlebodies. Finally, in Section~3 we describe how the $L$--moves adapt for deriving braid equivalences for virtual knots, welded knots, flat knots and for singular knots.

\smallbreak

More precisely: the diagrammatic equivalence in a specific topological setting, but also moves that are not allowed in the setting, define the corresponding braid isotopy, the way the closure of a braid is realized, and also the types of $L$--moves needed for the corresponding braid equivalence.

In the case of knots and links in knot complements, in c.c.o. $3$--manifolds and in handlebodies, the  basic idea was to represent the manifold by a braid in $S^3$ and then links and braids in the manifold by mixed links and mixed braids in $S^3$. Then to translate isotopy in these spaces in terms of mixed link isotopy in $S^3$, according to the hosting manifold each time. The main difference among the knot theories in these spaces is that in a handlebody a knot may not pass beyond the boundary of the handlebody from either end, and this is reflected both in the definition of the closure of a braid (see Figs.~27, 28) as well as in the corresponding braid equivalence. Further, in the case of c.c.o. $3$--manifolds we have an extra isotopy move (from the case of knot complements) coming from the handle sliding moves related to the surgery description of the manifold. All these are explained in Section~2.

Regarding the $L$--move braid equivalence corresponding to the knot isotopy in a given setting, there are two conceptual steps: first create the $L$--move braid equivalence on a geometric level. Then take into account the algebraic structures of the braids and the algebraic expressions of the $L$--moves (using the interpretation of an $L$--move as introducing an in-box crossing, see Fig.~6 for braids in $S^3$) and turn the geometric braid equivalence into an algebraic braid equivalence. The interpretation of an $L$--move as introducing an in-box crossing is also very crucial  in the search of the types of $L$--moves needed in a specific diagrammatic setting, as they are related to the types of kinks allowed in the given diagrammatic isotopy. For example, in the virtual domain we have $L$--moves introducing a real or a virtual crossing, facing to the right or to the left of the braid. Moreover, the presence of the two forbidden moves in the theory leads to the fact that the strands of an $L$--move cross the other strands of the braid only virtually, and also to a type of virtual $L$--move coming from a `trapped' virtual kink (see Definitions~10, 11 and~12). As another example, in the singular domain we do not have $L$--moves introducing a singular crossing, as the closure of such a move would contract to a kink with a singular crossing, and this is not an isotopy move in the theory. All these are discussed in Section~3.

Regarding the algebraic braid equivalences, surprisingly conjugation  by a braid generator may sometimes be achieved simply by $L$--moves within the braid. Then, such conjugations will not appear in the geometric statement of the braid equivalence.  But there are also situations where a braid generator, even if it has an inverse, cannot be conjugated. For example, the loop generators $a_i$ in a handlebody of genus $m\geq 2$ (see Subsection~2.4), which can be conjugated in the case of knot complements. Then, such conjugations will not appear even in the algebraic braid equivalence of the theory. Finally, there are braid structures where some
generators do not have inverses (e.g. the singular crossings $\tau_i$ in singular knot theory), so conjugation by these generators may not appear in the braid statement. Yet, a cutting line is permitted to cut a closed braid before or after a braid generator $x_i$. For
this reason we will adopt the following more general move, which is equivalent to conjugation if $x_i$ has an inverse.

\vspace{.05in}
\noindent $\bullet$ {\it Commuting in ${\mathcal B}_n$ :  $\alpha x_i \sim x_i \alpha$,}
\vspace{.05in}

\noindent where ${\mathcal B}_n$  denotes the set of braids with $n$ strands in the given diagrammatic knot theory. There may be also situations where the cutting line of the closure cannot cut anywhere
 (for example, singular knot theory in a handlebody). For this reason we will adopt the more general stabilization move:

\vspace{.05in}
\noindent $\bullet$ {\it Stabilization in $\cup_n{\mathcal B}_n$: $\alpha_1 \alpha_2 \sim \alpha_1 x_n^{\pm 1} \alpha_2$,}
\vspace{.05in}

\noindent  where $\alpha_1, \alpha_2 \in {\mathcal B}_n$ and $x_n \in {\mathcal B}_{n+1}$.
Stabilization by a real crossing may occur in all the above-mentioned knot theories. This is because the
isotopy move RI is permitted in all these theories. Yet, stabilization by a singular crossing is not
permitted in singular braids.

\smallbreak

The $L$--moves between braids are very fundamental and they provide a flexible conceptual center from which to deduce many results. They may be adapted to any diagrammatic setting where the notion of braid and diagrammatic isotopy is defined,  and they provide a uniform ground for formulating and proving braid equivalence theorems, first geometric and then algebraic.
 Moreover, the local algebraic versions of the braid equivalences in the various settings, as they arise using the $L$--moves, promise to be useful for constructing invariants of knots and links in these settings.

 As we said earlier, the differences in the various settings lie in the knot isotopies allowed or not allowed in the theory, and these will be reflected in the way the closure of a braid is defined and in the corresponding braid equivalences. For instance, one could work out braid equivalences for the virtual or the welded or the singular knot theory in knot complements, in c.c.o. $3$--manifolds or in handlebodies. In future work we intend to apply our methods to transverse knots and to knots in thickened surfaces.

 All the issues treated in this paper are expanded, discussed thoroughly and proved in detail in the author's monograph \cite{La4}.

\section{Braid equivalence in $S^3$}

A {\it link} of $k$ components is the homeomorphic image of $k$ copies of the circle into three-space and a knot is a link of one component. The central problem in classical Knot Theory is the complete classification of knots and links in $S^3$ up to ambient isotopy. Reidemeister \cite{Rd1} translated isotopy in space into a discrete diagrammatic equivalence relation on knot diagrams in the plane, generated by the three well--known Reidemeister moves RI, RII, RIII together with planar isotopy, see Fig.~1. Further, isotopy of oriented links uses all variations of the Reidemeister moves resulting in from all  different orientations of the arcs involved. In the sequel, whenever we say `knots' we shall mean both knots and links.

\begin{figure}
     \begin{center}
     \includegraphics[width=6cm]{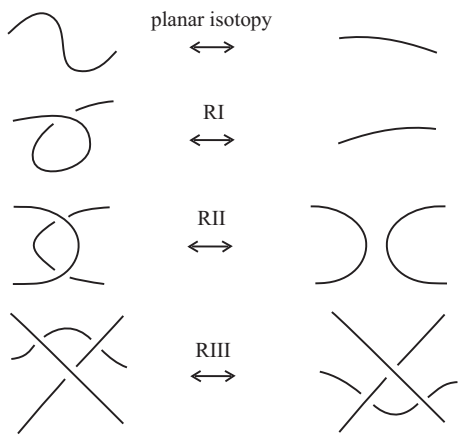}
     \caption{The diagrammatic isotopy moves }
     \label{}
\end{center}
\end{figure}

Braids, as topological objects, are similar to links, but they have an algebraic structure. A {\it braid} on $n$ strands  consists in $n$ arcs embedded in a thickened square $[0,1] \times [0,\epsilon] \times [0,1]$, such that the ends are arranged into $n$ collinear top endpoints in $[0,1] \times [0,\epsilon] \times \{1\}$ and  $n$ collinear bottom endpoints in $[0,1] \times [0,\epsilon] \times \{0\}$ and such that the embedding is monotonous with respect to the height function, that is, there are no local maxima or minima. Two braids are isotopic if one can be deformed to the other via an isotopy that preserves the braid structure. Braids shall be identified with braid diagrams (that is, regular projections on $[0,1] \times \{\epsilon\} \times [0,1]$) with no risk of confusion. Algebraically, the set of braids on $n$ strands forms the classical braid group $B_n$ \cite{Ar1,Ar2}, \cite{Ch}, whose generators $\sigma_1,\ldots, \sigma_{n-1}$ are the elementary positive crossings between consecutive strands and they satisfy the defining {\it braid relations}:
\begin{equation}\label{brels}
\sigma_{i}\sigma_{i+1}\sigma_{i}= \sigma_{i+1}\sigma_{i}\sigma_{i+1}\qquad \text{and}\qquad
\sigma_i\sigma_j = \sigma_j\sigma_i \quad \text{for}\quad \vert i-j\vert >1.
\end{equation}
The first relation is the most characteristic, and it reflects the braid isotopy move RIII (all arrows down).  Also, the fact that the generators $\sigma_i$ are invertible is reflected into the braid isotopy move RII. An excellent reference on braid groups is \cite{KT}.

By definition a braid has a natural orientation (from top to bottom). Taking the closure of a braid $\beta$, that is, joining its corresponding endpoints with simple unlinked arcs to the side of the braid yields an oriented link, denoted $\widehat{\beta}$. Conversely we have:

\begin{thm}
[J.~W.~Alexander, 1923] \label{alex}
 Any oriented link in $S^3$ is isotopic to the closure of a braid (not unique).
\end{thm}

Brunn \cite{Br} proved that any link has a projection with a single
multiple point; from which it follows immediately (by appropriate perturbations) that we can braid
any link diagram. Apart from Alexander \cite{Al}, other proofs of Theorem~\ref{alex}  have been given by  Morton \cite{Mo}, Yamada \cite{Ya}, Vogel \cite{Vo}, by the author with Rourke \cite{La1,LR1}, and by the author with Kauffman \cite{KL1}, where the emphasis is given on braiding virtual knots.

\smallbreak
We shall present here two braiding algorithms, that of \cite{La1,LR1} and that of \cite{KL1}, both conceptually very simple and both based on the same idea. The context of both braidings is that of \cite{La1,LR1}: Any  link
diagram can be arranged to be in general position with respect to the standard height function on
the plane. This means that it does not contain any horizontal arcs and it can be seen as a composition
of horizontal stripes, each containing either a local minimum or a local maximum or a crossing. The idea of the braiding  is on the one hand to keep the {\it down--arcs} of the diagram that are oriented downward and, on the other hand, to
eliminate the {\it up--arcs} that go upward and produce instead braid strands.
 First consider up--arcs that contain no crossings. Call such an arc in the diagram a {\it free up--arc}.

\begin{defn} \label{braiding} \rm
 A {\it braiding move} on a free up--arc is the following elimination operation of it:  We cut the arc at any point.  We then pull the two ends the upper upward and the lower downward, keeping them aligned, and so as to obtain a pair of corresponding braid strands, both running entirely {\it over\/}
the previously constructed tangle or entirely {\it under\/} it. View Fig.~2 for an abstract illustration. In Fig.~2 the up--arc is cut at the topmost point and the two ends are both pulled over the rest of the diagram, which is represented by the bold circle. The closure of the resulting tangle is a link diagram, obviously isotopic to the original one.
 Indeed, from the up--arc we created a stretched loop, which by move RI is isotopic to the original
free up--arc. Note that in Fig.~2 a faint arc is also illustrated inside the braid. It should be ignored in this definition; we shall refer to it later.
\end{defn}

\begin{figure}
     \begin{center}
     \includegraphics[width=12cm]{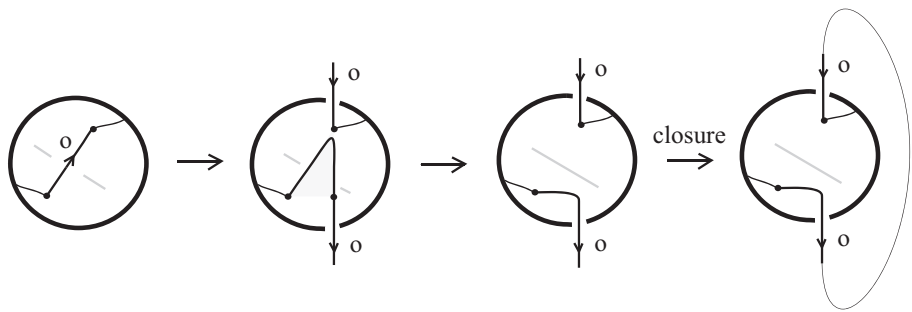}
     \caption{The  basic braiding move }
     \label{}
\end{center}
\end{figure}

In \cite{KL1}, before performing the braiding moves we prepare the diagram by rotating all crossings that contain up-arcs, so that the arcs that pass through any crossing are directed downward. There are two
types of rotation: If both arcs of the crossing go up, then we rotate the crossing by 180 degrees. If only one arc in the crossing goes up, then we rotate it by 90 degrees. See Figs.~3 and~4.
 These rotations may produce new free up--arcs. After adjusting all the crossings we then braid all the
free up--arcs. The resulting tangle is the desired braid, the closure of which is isotopic to
the original diagram.

\begin{figure}
     \begin{center}
     \includegraphics[width=5.3cm]{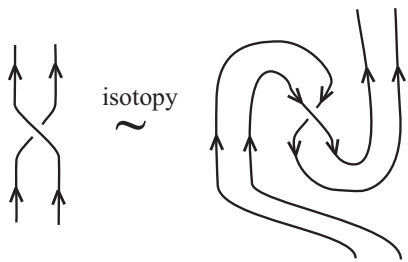}
     \caption{Full twist}
      \label{}
\end{center}
\end{figure}

\begin{figure}
     \begin{center}
     \includegraphics[width=3.7cm]{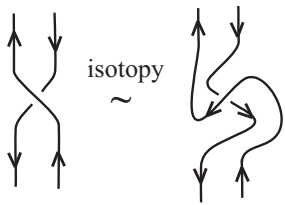}
     \caption{Half twist}
      \label{}
\end{center}
\end{figure}

In \cite{La1,LR1} we do not rotate the crossings. Instead, we apply the braiding move also to the up--arcs going through  crossings. The difference here is that we have to label each up--arc as {\it over} or {\it under} according to its position in the crossing, and then the resulting strands from the braiding move will be both over or under the rest of the diagram accordingly, see \cite{LR1} for details. For example, in Fig.~2 consider the up--arc creating a crossing with the faint arc. Then the up--arc is the over arc and it is labelled ``o", so the two new strands are stretched over the rest of the diagram.

\begin{rem} \rm
 The advantage of the braiding technique in \cite{KL1} is that it applies equally well to the virtual category, to flat virtuals, to welded knots and to singular knots, see \cite{KL1} for details. It applies, in fact, to all the categories in which braids are constructed. The braiding technique in \cite{La1,LR1} is preferred in situations where we need to analyze how the isotopy
moves on diagrams affect the final braids, see \cite{LR1,HL,KL1} for details.
\end{rem}

Trying now to use braids for studying knots we need to know when two braids correspond to isotopic knots. The Markov theorem gives an equivalence relation in the set of all braids $B_\infty := \cup_{n=1}^{\infty} B_n$ reflecting isotopy of knots:

\begin{thm}
[A. A. Markov, 1936]\label{markov}
 Two oriented links in $S^3$ are isotopic if and only if any two corresponding braids in $B_\infty$ differ by braid relations and the moves:

\noindent (i) Conjugation in $B_n$: \ $\sigma_i^{-1} \alpha \sigma_i \sim \alpha$,

\noindent (ii) Bottom stabilization in $B_\infty$: \ $\alpha \sim \alpha \sigma_n^{\pm 1}$, where $\alpha \in B_n$.
\end{thm}

 Conjugation can be generalized to the braid equivalence move {\it commuting} described in the Introduction. Indeed, conjugation is clearly related to isotopy move RII occurring in the
back side of a closed braid $\widehat{\alpha}$. Opening $\widehat{\alpha}$
with a cutting line, this line will either put both crossings of the move at the same side of the braid
(and then the initial RII move becomes the braid cancellation $\sigma_i^{-1}\sigma_i$) or it will
separate the two crossings, one at the top and one at the bottom of the braid, yielding a conjugate of $\alpha$.
Since the $\sigma_i$'s are invertible, commuting in this case is equivalent to conjugation.  Regarding the bottom stabilization, we just note that the more general stabilization move described in the Introduction together with move (i) implies bottom stabilization. Stabilization is related to the isotopy move RI, since the closure of a braid with a stabilization move is a knot with a kink at that place.

The Markov theorem is not easy to prove. Its proof is very technical and it depends on the braiding algorithm used. Markov \cite{Ma} used Alexander's braiding algorithm. The Markov theorem was originally stated by Markov \cite{Ma} with {\it three} braid moves and then Weinberg \cite{We} reduced them to the {\it two} braid moves of Theorem~\ref{markov}.   Birman \cite{Bi1} gave a more rigorous proof, filling in all details of Markov's original proof, using a more rigorous version of Alexander's braiding algorithm. Bennequin \cite{Ben} gave another proof using contact topology. After the discovery of the Jones polynomial \cite{Jo} --the first knot invariant that was constructed using braids and the Markov theorem-- other proofs of Theorem~\ref{markov} have been given by: Morton \cite{Mo} using his threading algorithm, by the author with Rourke \cite{La1,LR1} using the braiding algorithm of \cite{La1,LR1} and the new type of braid move, the $L$--move (see below), by Traczyk \cite{Tr} using Vogel's braiding algorithm, and by Birman and Menasco using Bennequin's ideas.

\smallbreak
 In \cite{LR1} using the $L$--moves we proved, in fact, the following sharper version of Theorem~\ref{markov} (cf. \cite{LR1}, Thm. 2.3).

\begin{thm}[One-move Markov theorem] \label{lthm}
 Two oriented links in $S^3$ are isotopic if and only if any two corresponding braids differ by braid relations
 and the $L$-moves.
\end{thm}

\begin{defn}\label{lmove} \rm An {\it $L$--move} consists in cutting a strand of the braid at a point, then pulling the two ends of the cut, the upper downward and the lower upward, both running entirely {\it over\/} the braid or entirely {\it under\/} it, and keeping them aligned (on the vertical line of the cutpoint) so as to obtain a  new pair of corresponding braid strands.  Thus there are two types of $L$--moves, an {\it under $L$--move}, or {\it $L_u$--move}, and an {\it over $L$--move}, or {\it $L_o$--move}. Fig.~5 illustrates an example of the two $L$--moves.
 Further, by a small braid isotopy that does not  change the relative positions of endpoints, an $L$--move is equivalent to adding an in-box crossing (positive or negative) facing the right-hand side or the left-hand side of the braid box.
\end{defn}

So, the $L$--move generalizes the stabilization move. View last instance of Fig.~5, where a positive in-box crossing is formed and two dots are placed to indicate the starting points of the $L_o$--move. See also Fig.~6 for an abstract illustration. Closing the new pair of strands created from an $L$--move results in a kink in the tangle diagram, so it corresponds to the isotopy move RI.

\begin{figure}
     \begin{center}
     \includegraphics[width=15cm]{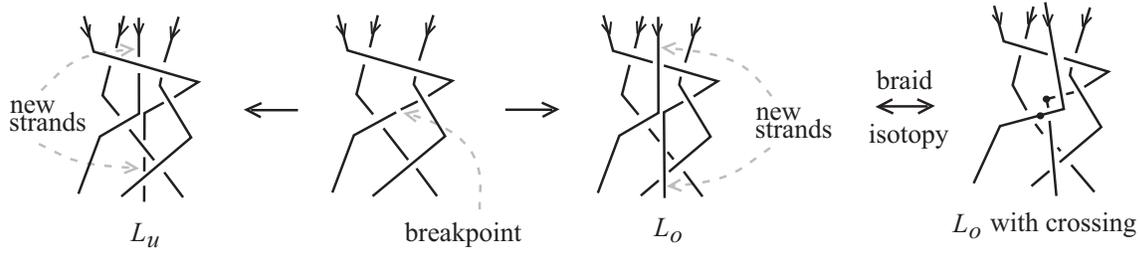}
     \caption{The $L$--moves }
     \label{}
\end{center}
\end{figure}

\begin{rem} \rm The version of an $L$--move introducing an  in-box crossing gives rise to the following algebraic expression for an $L_o$--move and an $L_u$--move respectively:
\begin{equation}\label{alglo}
\alpha=\alpha_1\alpha_2 \stackrel{L_o}{\sim}
\sigma_i^{-1}\ldots \sigma_n^{-1} \alpha_1' \sigma_{i-1}^{-1}\ldots
\sigma_{n-1}^{-1}\sigma_n^{\pm 1}\sigma_{n-1} \ldots \sigma_i
\alpha_2' \sigma_n \ldots \sigma_i
\end{equation}
\begin{equation}\label{alglu}
\alpha=\alpha_1\alpha_2 \stackrel{L_u}{\sim}
\sigma_i\ldots \sigma_n \alpha_1' \sigma_{i-1}\ldots
\sigma_{n-1}\sigma_n^{\pm 1}\sigma_{n-1}^{-1}\ldots\sigma_i^{-1}
\alpha_2' \sigma_n^{-1}\ldots\sigma_i^{-1}
\end{equation}

\noindent  where $\alpha_1$,
 $\alpha_2$ are elements of $B_n$ and $\alpha_1'$,
 $\alpha_2' \in B_{n+1}$ are obtained from $\alpha_1$,
 $\alpha_2$ by replacing each $\sigma_j$ by $\sigma_{j+1}$ for
 $j=i,\ldots,n-1$.
 \end{rem}

\begin{figure}
     \begin{center}
     \includegraphics[width=14cm]{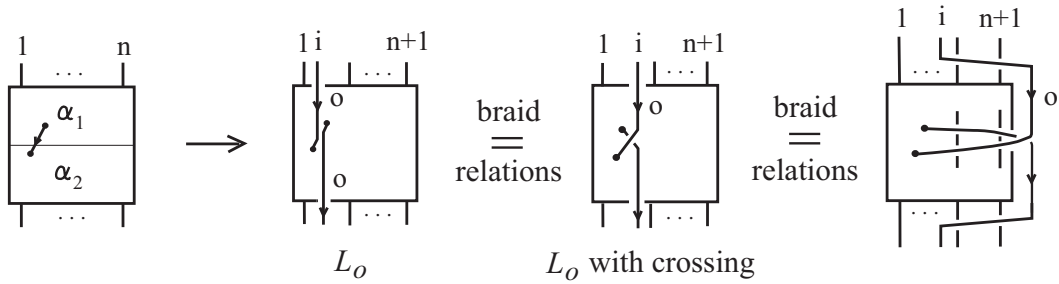}
     \caption{ $L_o$--move, $L_o$--move with crossing and its algebraic expression }
     \label{}
\end{center}
\end{figure}

Theorems~\ref{markov} and~\ref{lthm} imply that braid conjugation follows from the $L$--moves. This fact is also directly illustrated in Fig.~7.

\begin{figure}
     \begin{center}
     \includegraphics[width=14cm]{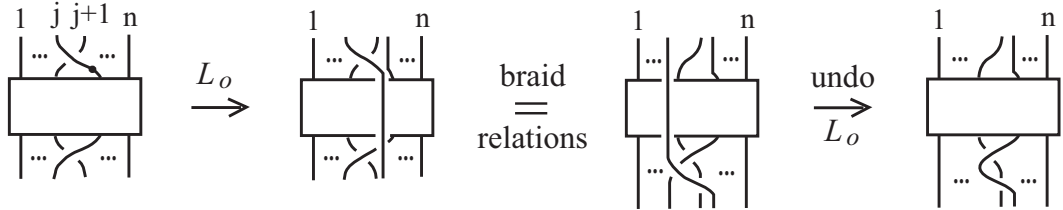}
     \caption{ Braid conjugation follows from the $L$--moves }
     \label{}
\end{center}
\end{figure}

\begin{rem} \rm
 Theorem~\ref{lthm} has three main advantages: First, there is only one type of braid equivalence move. Second, the $L$--moves are very local and fundamental. Thus, adapted each time to a given diagrammatic setting (where diagrammatic isotopy is given and the notion of a braid can be defined), they may provide a uniform ground for formulating and proving geometric analogues of the Markov theorem.  Third, due to the fact that the $L$--moves have algebraic expressions, these geometric statements can be turned into algebraic statements,
according to the braid structures of the given situation (which, in turn, depend on the diagrammatic equivalence and the manifold where the diagrams live).
\end{rem}

\section{Braid equivalence in  knot complements, c.c.o. 3-manifolds and handlebodies}

It is a well--known result in Topology
that any c.c.o. $3$-manifold can be obtained via surgery on an unoriented framed link in $S^3$ with integral framings. The intermediate stage of surgery is the construction of the complement of the surgery link. By the term `knot complement' we shall be referring throughout to both knot and link complements. Another category of bounded 3-dimensional manifolds that give rise to c.c.o. 3-manifolds are the handlebodies, via the Heegaard decomposition. The special case of the solid torus is the only manifold common in both categories, and its knot theory has been studied quite extensively from various viewpoints.

\subsection{The common diagrammatic setting}

 As we shall establish in this section, knots and braids in all three categories of 3-dimensional manifolds  may be studied via very similar diagrammatic representations in $S^3$ and very similar techniques. The differences lie in the knot isotopies allowed or not allowed in each manifold, and these will be reflected in the way we define the closure of a braid and in the corresponding braid equivalences. The underlying braid moves that lead to algebraic braid equivalences in all these manifolds are the $L$--moves.

\smallbreak
 Let $S^3 \backslash K$ be the complement of the oriented link $K$ in $S^3$. By the Alexander theorem and the
definition of ambient isotopy, $S^3 \backslash K$ is homeomorphic to $S^3 \backslash \widehat{B}$,  where $\widehat{B}$ is the closure of some braid $B$  and it  is isotopic to $K$. Throughout this section the braid $B$ will remain fixed. Let, further, $L$ be an oriented link in $S^3 \backslash \widehat{B}$. Fixing $\widehat{B}$ pointwise on its projection plane we may represent $L$ unambiguously  by the {\it mixed link} $\widehat{B}\bigcup L$ in $S^3$, that is,
a link in $S^3$  consisting of the {\it fixed part} $\widehat{B}$ and the {\it moving part} $L$
that links with $\widehat{B}$, view Fig.~8a for an example.  A {\it mixed link diagram} is a diagram of $\widehat{B}\bigcup L$ projected on the plane of $\widehat{B}$ which is equipped with the top-to-bottom direction of $B$.

\begin{figure}
     \begin{center}
     \includegraphics[width=15cm]{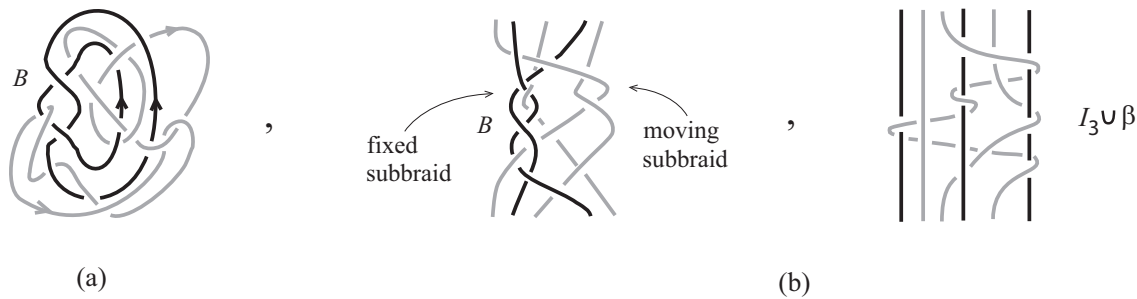}
     \caption{A mixed link and two mixed braids related to knot complements and c.c.o. 3-manifolds}
     \label{}
\end{center}
\end{figure}

\smallbreak

Let now $M$ be a c.c.o. $3$-manifold obtained by surgery on an unoriented framed link in $S^3$. Without loss of generality this link is the closure $\widehat B$ of some braid $B$.  We shall refer to
$\widehat B$ as the {\it surgery link} and we shall be writing $M=\chi (S^3, \widehat B)$. So $M$ may be represented in $S^3$ by the framed braid $\widehat B$. Further, fixing $\widehat B$ pointwise, we have that links in $M$ can be unambiguously represented by mixed links in $S^3$ --exactly as in the case of knot complements.
So, we follow the same diagrammatic setting.

\smallbreak
Trying to find an analogue of the Markov theorem for links in $S^3 \backslash \widehat{B}$ or $\chi (S^3, \widehat B)$ we observe first that  $\widehat{B}$ must remain fixed. So, we need a braiding process for mixed link diagrams that maps $\widehat{B}$ to $B$ and not to an $L$--equivalent braid or even to a conjugate of $B$. Such braids are called geometric mixed braids. A {\it geometric mixed braid} is a braid $B\bigcup \beta$ on $m+n$ strands that contains the subbraid $B$ as a fixed subbraid. View Fig.~8b for two examples. The subbraid $B$ represents the manifold where the knot theory lives in and it shall be called {\it the fixed subbraid} in the case of knot complements and {\it the surgery braid} in the case of c.c.o. 3-manifolds, while the subbraid $\beta$ that complements $B$ shall be called {\it the moving subbraid} and it represents the link in the manifold. These braids were introduced in \cite{La1,LR1}. The closure of a geometric mixed braid is defined to be the standard closure of the mixed braid seen as a classical braid on $m+n$ strands. As we shall see below, every mixed link may be braided to a geometric mixed braid.

\smallbreak

Let $H_m$ denote the handlebody of genus $m$. The {\it handlebody} of genus $m$ is usually defined as $(\mbox{\it a closed disc} \setminus m \, \mbox{\it open discs})\times I$, where $I$ is the
unit interval. Equivalently, $H_m$ can be defined as $(S^3 \setminus \mbox{\it an open tubular
neighbourhood of} \ I_m)$, where $I_m$ denotes here the identity braid on $m$
infinitely extended strands, all meeting at the point at infinity, see Fig.~9a. Thus $H_m$
may be represented in $S^3$ by the usual identity braid  $I_m$. Let now $L$ be an oriented link in $H_m$. Then, fixing $I_m$ pointwise on its projection plane, $L$ may be represented unambiguously  by the mixed $(m,m)$-tangle $I_m \bigcup L$ in $S^3$ which, by abuse of language, we shall call {\it mixed link} (view Fig.~9b).  The subbraid $I_m$ shall be called  {\it the fixed part}
and $L$  {\it the moving part} of the mixed link. A {\it mixed link diagram}
is then a diagram of $I_m \bigcup L$ projected on the plane of $I_m$, which is equipped with the
top-to-bottom direction.

\begin{figure}
     \begin{center}
     \includegraphics[width=12cm]{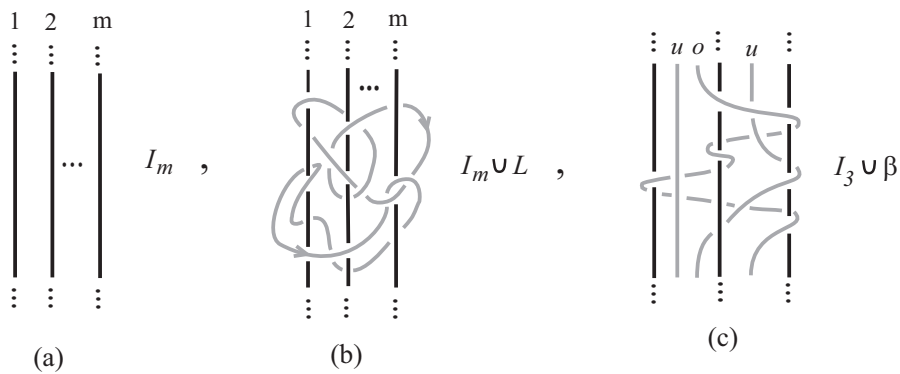}
     \caption{ Representing $H_m$ - a mixed link and a mixed braid related to a handlebody }
     \label{}
\end{center}
\end{figure}

Now, applying the braiding algorithm of \cite{La1,LR1}, by which all down--arcs remain fixed, we can easily braid an oriented mixed link $I_m \bigcup L$ in $S^3$ to a geometric mixed braid $I_m \bigcup \beta$, as defined above. In the case of a handlebody, though, the strands of the fixed subbraid $I_m$ are thought infinitely extended. This forces the strands of the moving subbraid to be labelled with labels `u' and `o' for `under' and `over' (and this is related to the definition of the closure, see details in the subsection on handlebodies). These braids were introduced in \cite{HL}. View Fig.~9c for an example.

\smallbreak
 We note that, in all settings, if we remove the fixed braid from a mixed link or a geometric mixed braid we are left with an oriented link resp. a braid in $S^3$.

\subsection{The common algebraic setting}

An {\it algebraic mixed braid} on $n$ moving strands is a braid on $m+n$ strands such that the first $m$ strands form the identity braid $I_m$. View Fig.~10 for an example. These braids were introduced in \cite{La2}. We denote the set of all such braids $B_{m,n}$. Clearly, if in the definition of a geometric mixed braid $B$ stands for the identity braid $I_m$, then the set of algebraic mixed braids $B_{m,n}$ is a subset of the corresponding set of geometric mixed braids.

\begin{figure}
\begin{center}
\includegraphics[width=3.2cm]{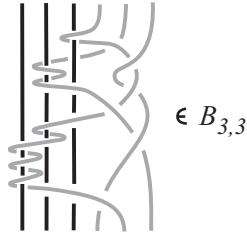}
\end{center}
\caption{ An algebraic mixed braid }
\label{}
\end{figure}

\smallbreak

The set $B_{m,n}$ is a subset of $B_{m+n}$ closed under concatenation and under inverses, so it has a  group structure and shall be called {\it the mixed braid group}. The mixed braid group is generated by the classical crossings $\sigma_j$ and by the `loops' $a_i$ as illustrated in Fig.~11. Moreover, $B_{m,n}$ has the following presentation (cf. \cite{La2}):
\begin{equation}
B_{m,n} = \left< \begin{array}{ll}  \begin{array}{l}
a_1, \ldots, a_m,  \\
\sigma_1, \ldots ,\sigma_{n-1}  \\
\end{array} &
\left|
\begin{array}{l} \sigma_k \sigma_j=\sigma_j \sigma_k, \ \ |k-j|>1   \\
\sigma_k \sigma_{k+1} \sigma_k = \sigma_{k+1} \sigma_k \sigma_{k+1}, \ \  1 \leq k \leq n-1  \\
{a_i} \sigma_k = \sigma_k {a_i}, \ \ k \geq 2, \   1 \leq i \leq m    \\
 {a_i} \sigma_1 {a_i} \sigma_1 = \sigma_1 {a_i} \sigma_1 {a_i}, \ \ 1 \leq i \leq m  \\
 {a_i} (\sigma_1 {a_r} {\sigma^{-1}_1}) =  (\sigma_1 {a_r} {\sigma^{-1}_1})  {a_i}, \ \ r < i.
\end{array} \right.  \end{array} \right>
\end{equation}

\begin{figure}
\begin{center}
\includegraphics[width=4.9in]{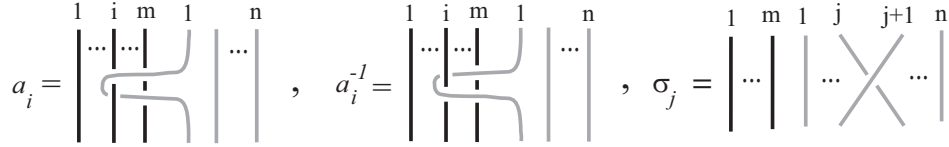}
\end{center}
\caption{ The `loops' $a_i, \ {a^{-1}_i}$ and the crossings $\sigma_j$ }
\label{}
\end{figure}

\smallbreak

As it turns out \cite{La2,LR1,HL} the mixed braid group $B_{m,n}$ is the natural algebraic companion to oriented knots and links in knot complements, c.c.o. 3-manifolds and handlebodies. More precisely, the groups $B_{m,n}$ are the appropriate braid structures for studying knots and links in the complement
of the $m$--unlink or a connected sum of $m$ lens spaces or in a handlebody of genus
$m$. For arbitrary knot complements and c.c.o. 3-manifolds we have specific cosets
$C_{m,n}$ of these groups in the groups $B_{m+n}$, which are determined by the specific hosting manifold, represented in $S^3$ by the braid $B$  \cite{La2}. Then an element of the related coset $C_{m,n}$ consists in an elements of the group $B_{m,n}$ followed by the natural embedding of $B$. Yet, as we shall explain below, by parting  and combing mixed braids the braid equivalences in all these spaces can be formulated only in terms of the mixed braid groups $B_{m,n}$.
 The capstone for studying braid equivalence in all above-mentioned spaces is the following result (\cite{LR1}, Thm. 4.7).

 \begin{thm}[Relative version of the $L$--move braid equivalence]\label{relv} Let $L_1$, $L_2$ be oriented link diagrams which both contain a common braided portion
$B$. Suppose that there is an isotopy of $L_1$ to $L_2$ which finishes with a homeomorphism fixed
on $B$. Suppose further that $B_1$ and $B_2$ are braids obtained from our braiding process applied
to $L_1$ and $L_2$ respectively. Then $B_1$ and $B_2$ are equivalent by braid isotopy and by by $L$--moves that do not
affect the common braided portion $B$.
\end{thm}

 Using Theorem~\ref{relv}, geometric versions of the Markov theorem were proved for links in knot complements and for links in arbitrary c.c.o. $3$-manifolds by means of $L$--move analogues (cf. \cite{LR1}, Thms. 5.5 and 5.10). These results were consequently turned into algebraic statements in \cite{LR2} (Thms.~4 and~5), via the braid groups $B_{m,n}$. Theorem~\ref{relv} led also to an $L$--move analogue of the Markov theorem for knots in a handlebody as well as to algebraic analogues of it using the braid groups $B_{m,n}$ (cf. \cite{HL}, Thms.~3, 4 and 5). We shall now explain all this in more detail.

\subsection{Braid equivalence in knot complements}

Two oriented links $L_1$, $L_2$ are isotopic in $S^3 \backslash \widehat{B}$ if and only if the mixed links $\widehat{B}\bigcup L_1$ and $\widehat{B}\bigcup L_2$ are isotopic in $S^3$ by an ambient
isotopy which keeps $\widehat{B}$ pointwise fixed. In terms of diagrams, mixed link isotopy is generated by the ordinary Reidemeister moves for the moving part and by the mixed moves RII, RIII of Fig.~12 involving also the fixed part (see \cite{LR1} for details).

\begin{figure}
     \begin{center}
     \includegraphics[width=9cm]{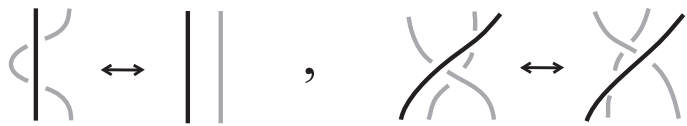}
     \caption{ The mixed Reidemeister moves}
     \label{}
\end{center}
\end{figure}

Further, in \cite{La1,LR1} we proved the following (cf. Theorem~5.3 \cite{LR1}):

\begin{thm}[Braiding theorem for $S^3 \backslash \widehat B$ and $\chi (S^3, \widehat B)$]\label{braidcompl} Any oriented link $L$ in $S^3 \backslash \widehat B$ can be represented in $S^3$ by some geometric mixed braid $B\bigcup \beta$, the closure of which is isotopic in $S^3$ to the mixed link  $\widehat B \bigcup L$.
 \end{thm}

\begin{proof} The fixed part  $\widehat{B}$ may be viewed in $S^3$ as the braid $B$ union an arc, $k$ say, at infinity. This arc is the identification of the two horizontal arcs containing the endpoints
of $B$ and thus it realizes the closure of $B$, view Fig.~13. In Fig.~13a the braid $B$ is drawn curved but
this is just for the purpose of picturing it. Let $L$ be a link in $S^3 \backslash \widehat B$. By general position $L$ misses a small regular neighbourhood of $k$, $N(k)$, and
therefore it can be isotoped into the complement of $N(k)$ in $S^3$. By expanding $N(k)$ we can
view its complement as the cylinder $T=D^2\times I$ that contains $B$ (Fig.~13b). We then
apply the braiding of \cite{LR1} in $T$. This will leave $B$ untouched but it will braid $L$ in $T$.
Finally we cut $\widehat{B}$ open along $k$ in order to obtain a geometric mixed braid.
\end{proof}

Let $B \bigcup \beta$ be a braid in $S^3 \backslash \widehat B$. A {\it geometric $L$--move} consists in cutting a strand of the subbraid $\beta$ at some point and then performing an $L$--move as described in Definition~\ref{lmove}.
Using, now, Theorem~\ref{relv} we have the following  $L$--move version of the Markov theorem for knot complements (cf. Theorem~5.5 \cite{LR1}).

\begin{thm} [Geometric mixed braid equivalence for $S^3 \backslash \widehat B$]\label{lgeom} Two  oriented links $L_1, L_2$ in $S^3 \backslash \widehat B$ are isotopic  in $S^3 \backslash \widehat B$  if and
only if any two corresponding mixed braids $B \bigcup \beta_1$ and $B\bigcup \beta_2$ in $S^3$ differ by braid isotopies and by $L$--moves that do not touch the fixed subbraid $B$.
 \end{thm}

The passage from Theorem~\ref{relv} to Theorem~\ref{lgeom} is by showing that, if during the isotopy the hypothetical closing arc $k$ is crossed (see Fig.~13c), then the corresponding mixed braids differ by $L$--moves applied only on the moving subbraids.

\begin{figure}
     \begin{center}
     \includegraphics[width=12.5cm]{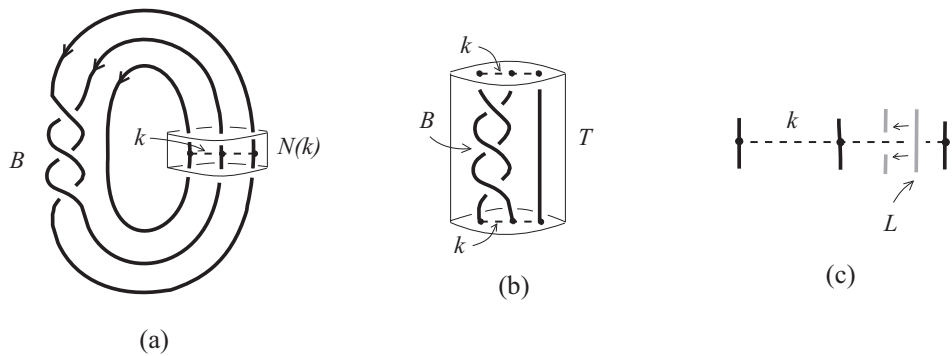}
     \caption{ Crossing the closing arc $k$}
     \label{}
\end{center}
\end{figure}

We would like now to pass from the geometric braid equivalence to an algebraic statement. Our strategy is the following. We first part the mixed braids and we translate the  $L$--equivalence of Theorem~\ref{lgeom} to an equivalence of parted mixed braids. {\it Parting} a geometric mixed braid $B \bigcup\beta$ on $m+n$ strands means to separate its endpoints into two different sets, the first $m$ belonging to the subbraid $B$ and the last $n$ to $\beta$, and so that the resulting braids  have isotopic closures.
 The generators of the groups $B_{m,n}$ become already apparent in the parted mixed braid equivalence. We then comb the parted mixed braids and we translate the parted mixed braid equivalence to an equivalence of algebraic mixed braids.  {\it Combing} a parted mixed braid means to separate the fixed subbraid from the moving subbraid using mixed braid isotopy.

\begin{lem}\label{parting} Every mixed braid may be represented by a parted mixed braid with
isotopic closure (cf. Lemma 1 \cite{LR2}).
\end{lem}

\begin{proof}  Let  $B \bigcup\beta$ be a geometric mixed braid. We attach arbitrarily arrays of labels `o' or `u' to each pair of corresponding strands of the moving subbraid $\beta$, with as many entries as the number of fixed strands on their right. We then pull each pair of corresponding moving strands to the right, {\it over\/} or {\it under\/} each  strand of $B$ lying on their right, according to the  label in the
array of the pair. We start from the rightmost pair, respecting the position of the
endpoints. See the first two illustrations of Fig.~14 for a parting of an abstract geometric mixed
braid. Obviously, the closures of the initial and of the parted mixed braids are isotopic (they differ by
planar isotopy and by mixed RII moves).
\end{proof}

Fig.~14 illustrates two different partings of an abstract geometric mixed braid, an arbitrary one and the standard parting, whereby all strands are pulled {\it over\/}. We shall use the same notation for parted mixed braids as for geometric mixed braids.

\begin{figure}
\begin{center}
\includegraphics[width=5.5in]{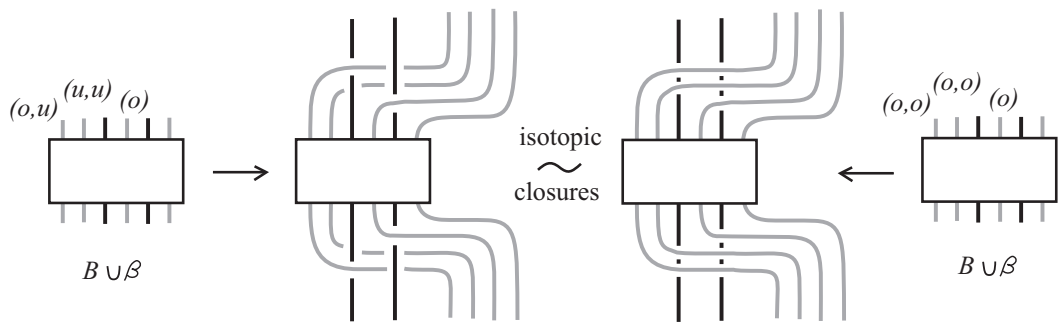}
\end{center}
\caption{ Parting a mixed braid -- the standard parting }
\label{}
\end{figure}

\smallbreak
Let now $C_{m,n}$ denote the set of parted mixed braids on $n$ moving strands, related to $S^3 \backslash \widehat B$. Note that for different fixed braids $B$ the corresponding sets of parted mixed braids will be also different. Yet, we will use the same notation for them with no risk of confusion. By adding an extra moving strand on the right of a parted mixed braid on  $n$ moving
strands, $C_{m,n}$ embeds naturally into $C_{m,n+1}$.  Let
$C_{m,\infty} := \bigcup_{n=1}^{\infty} C_{m,n}$ denote the disjoint union of  all sets
$C_{m,n}$. We define below some moves in  $C_{m,\infty}$.

\begin{defn}{\rm
{\it (1) \ Loop conjugation\/} of a parted mixed braid in $C_{m,n}$ is
its  concatenation from above by a loop $a_i$  (or by ${a^{-1}_i}$) and from
below by  ${a^{-1}_i}$ (corr.  $a_i$).

\vspace{.05in}
\noindent {\it (2) \  Markov conjugation\/} of a parted mixed braid  in $C_{m,n}$ is
its concatenation from above  by a crossings $\sigma_j$ (or by ${\sigma^{-1}_j}$) and from below by
${\sigma^{-1}_j}$  (corr.  $\sigma_j$).

\vspace{.05in}
\noindent {\it (3)} \ A {\it parted $L$--move\/} is defined to be  an  $L$--move  between
parted mixed braids. As for links in $S^3$, a parted $L$--move may be considered to introduce an in--box crossing (view left-hand illustration of Fig.~16).

\vspace{.05in}
\noindent  {\it (4)} \ A  {\it stabilization move} is the insertion of an extra strand with a crossing
${\sigma^{\pm 1}_n}$  at the right hand side of a parted mixed braid on $n$
moving strands or the reverse operation. View Fig.~15.

 The relations among moves (2), (3) and (4) are analogous to links in $S^3$ and are illustrated in Figs.~16 and 17.}
\end{defn}

\begin{figure}
\begin{center}
\includegraphics[width=3.2in]{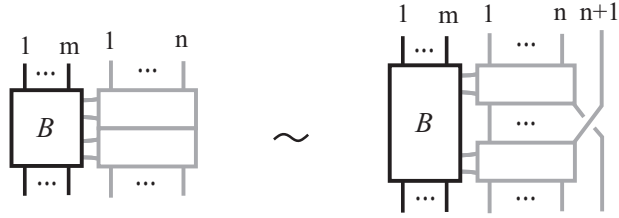}
\end{center}
\caption{ The stabilization move }
\label{}
\end{figure}

\begin{figure}
\begin{center}
\includegraphics[width=3.2in]{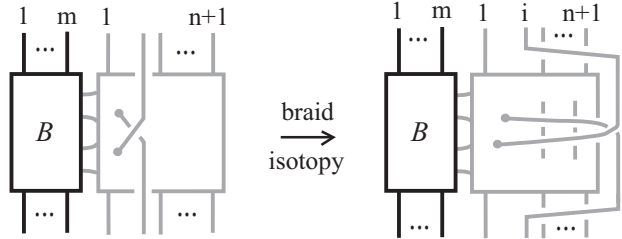}
\end{center}
\caption{ Parted $L$--move $\longleftrightarrow$ stabilization move and Markov conjugation }
\label{}
\end{figure}

\begin{figure}
\begin{center}
\includegraphics[width=5.7in]{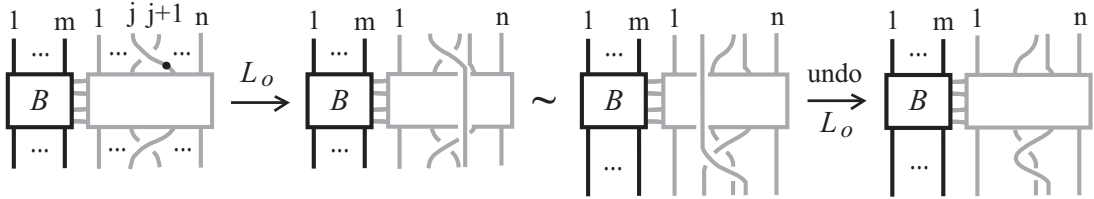}
\end{center}
\caption{ Markov conjugation is a composition of $L$--moves }
\label{}
\end{figure}

The lemma below says that, up to loop conjugations, we can choose any parting for a geometric mixed braid. Indeed, Fig.~18 illustrates that the change of only one label of one moving strand is realized by a loop conjugation.

\begin{lem} Two partings of a geometric mixed braid differ by loop conjugations (cf. Lemma~2 \cite{LR2}).
\end{lem}

\begin{figure}
\begin{center}
\includegraphics[width=5.7in]{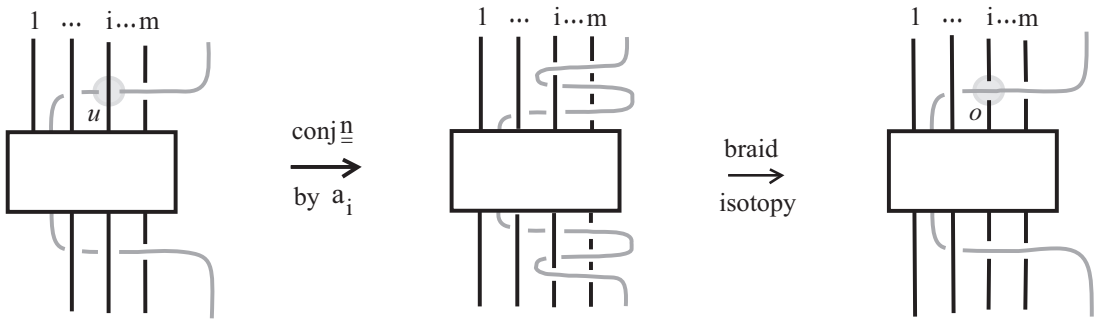}
\end{center}
\caption{ Change of parting labels $\longleftrightarrow$ conjugation by $a_i$ }
\label{}
\end{figure}

Finally, we need the following lemma:

\begin{lem} A  mixed braid with  an $L$--move performed can be parted to a parted mixed
braid with a parted $L$--move performed (cf. Lemma 3 \cite{LR2}.)
\end{lem}

\begin{proof}
 If the $L$--move is an $L_o$--move we choose to part its strands by pulling them to the right and over all other strands in between.  Then the crossing of the $L$--move slides over to the right by a braid isotopy. See Fig.~19. If the $L$--move is an $L_u$--move we pull the two strands under the fixed strands in between.
\end{proof}

\begin{figure}
\begin{center}
\includegraphics[width=4in]{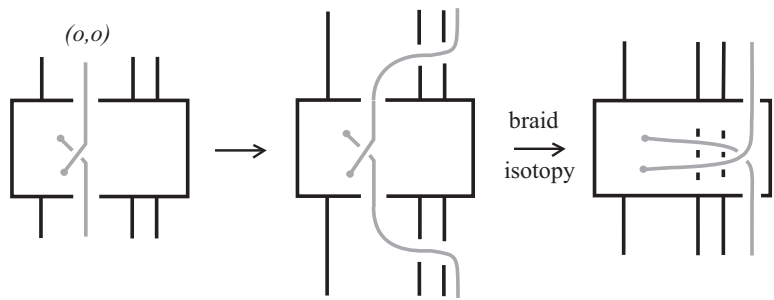}
\end{center}
\caption{   An $L_o$--move + parting = a parted $L_o$--move }
\label{}
\end{figure}

We are now in a position to state two versions of the braid equivalence  in $S^3 \backslash \widehat B$ for
parted mixed braids.

\begin{thm}[Parted mixed braid equivalence for $S^3 \backslash \widehat B$]\label{partedcompl} Two  oriented links in $S^3 \backslash \widehat B$ are isotopic if and only if
any two corresponding parted mixed braids in $C_{m,\infty}$ differ by  a finite sequence of braid isotopies, parted
$L$--moves and loop conjugations.

Equivalently,  by a  finite
sequence of braid isotopies, stabilization moves, Markov conjugations and loop conjugations.
\end{thm}

\begin{proof}
Clear by Lemmas~2 and 3 and by Figs.~16 and 17.
\end{proof}

We shall now move to a completely algebraic statement for the braid equivalence in $S^3 \backslash \widehat B$. Unless the fixed subbraid $B$ represents the complement of the $m$--unlink or a connected sum of $m$ lens spaces or a handlebody of genus $m$, $B$ is not the identity braid on $m$ strands and, so, concatenating two
elements of $C_{m,n}$ is not a closed operation, since it alters the braid description of
the manifold.  So, the set $C_{m,n}$ of parted mixed braids is not a subgroup of
$B_{m+n}$. Yet, we have the following result.

\begin{lem} For the fixed subbraid $B$ on $m$ strands the set $C_{m,n}$ is a coset of $B_{m,n}$ in $B_{m+n}$ (cf. Proposition~1 \cite{La2}.)
\end{lem}

\begin{proof} Let  $A\in C_{m,n}$. We shall show that $A$ can be written as a product $\alpha\, B$, where  $\alpha$ an algebraic mixed braid in $B_{m,n}$  followed  by the fixed subbraid $B$ embedded naturally in $B_{m+n}$. Indeed, we notice first that, by symmetry, Artin's combing for pure braids can be also applied starting from the bottom of the braid: first the pure braiding of the 2nd strand with the 1st and then remaining fixed, then the pure braiding of the 3rd strand with the 1st and the 2nd one and then remaining fixed, and so on.
 So, we multiply $A$ from the top with a braid $p$ on the last $n$ strands and with a
braid $\sigma$ on the first $m$ strands, such that $pA\sigma$ is a pure braid in  $B_{m+n}$. Then
we apply to it Artin's combing as above. This will separate $A$ into two parts: the top one (multiplied from the top with $p^{-1}$), an element of  $B_{m,n}$, representing the link in $S^3 \backslash \widehat B$, followed by the fixed braid $B$  embedded in $B_{m+n}$ (we have multiplied from the bottom with $\sigma^{-1}$). Thus, the set  $C_{m,n}$ of combed mixed braids is a coset of $B_{m,n}$ in $B_{m+n}$.
 \end{proof}

 A parted mixed braid in the form $\alpha\, B$, where $\alpha\in B_{m,n}$, shall be called a {\it combed mixed braid}. View Fig.~20 for an abstract illustration of a combed mixed braid.

\bigbreak

\begin{figure}
\begin{center}
\includegraphics[width=5cm]{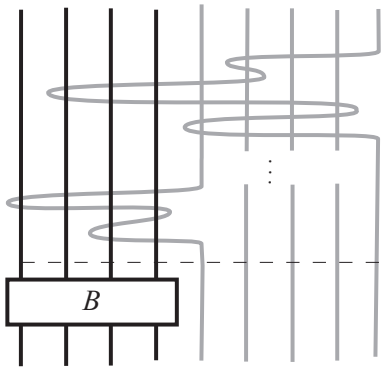}
\end{center}
\caption{ A combed mixed braid }
\label{figure18}
\end{figure}

In view of the above result we would like to restate the braid equivalence of Theorem~\ref{partedcompl} for parted mixed braids into an equivalence relation between their corresponding algebraic representatives after
combing. For this we need to understand how exactly the combing is done and how it affects
the parted braid equivalence moves.
Indeed, if we regard  a parted mixed braid as an element of the classical braid group $B_{m+n}$ and  if $\Sigma_k$ denotes the crossing between the $k$th and the $(k+1)$st strand of the fixed subbraid, then the crossings $\sigma_j$ of the moving part commute with the crossings of the fixed part, so they are not affected by the combing. Thus the only generating elements
of the moving part that are affected by the combing are the loops $a_i$.  Fig.~21 illustrates the relations between $\Sigma_1$ and the loops $a_i$. For  arbitrary $\Sigma_k$ the relations are completely analogous. Cf. \cite{LR2} for details.

\begin{figure}
\begin{center}
\includegraphics[width=5in]{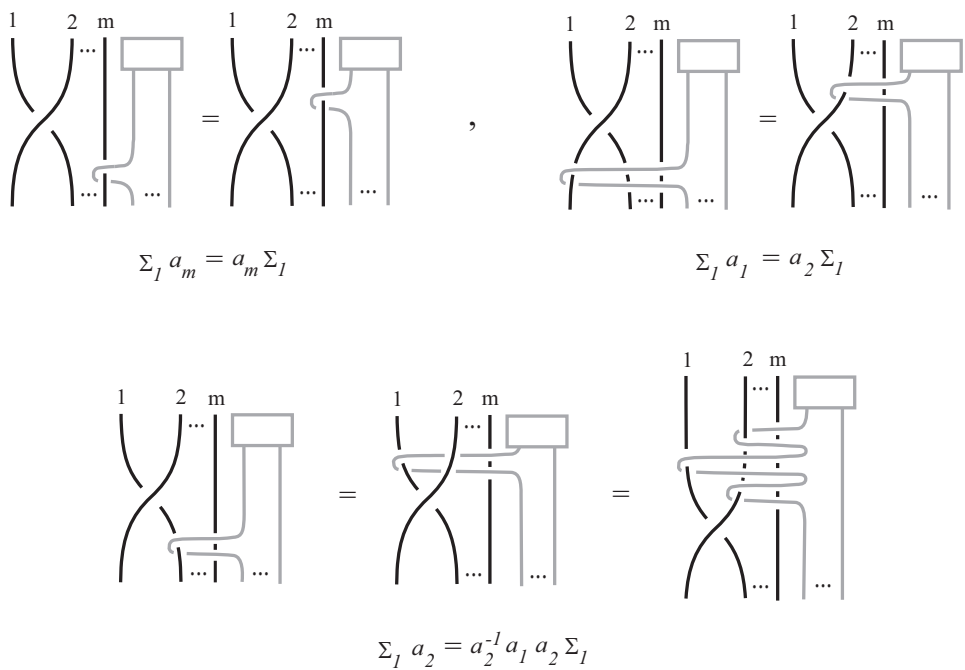}
\end{center}
\caption{ Combing the $a_i$'s to the top }
\label{figure19}
\end{figure}

\smallbreak

The group $B_{m,n}$ of algebraic mixed braids embedds naturally into the group $B_{m,n+1}$ and we shall denote  $B_{m,\infty} := \bigcup_{n=1}^{\infty} B_{m,n}$ the disjoint union of  all algebraic mixed braid groups. In $B_{m,\infty}$ we define the following moves.

\begin{defn}\label{algmoves} \rm
\noindent {\it (1) \  Algebraic Markov conjugation\/} is a Markov conjugation between algebraic mixed braids, and it has the algebraic expression:
\[
\alpha \sim {\sigma^{\pm 1}_j} \alpha {\sigma^{\mp 1}_j}
\]
 where $\alpha, \sigma_j \in B_{m,n}$.


\noindent  {\it (2)} \ An {\it algebraic stabilization move\/} is a stabilization move between algebraic mixed braids, and it has the algebraic expression:
\[
{\alpha}_1 {\alpha}_2 \sim {\alpha}_1{\sigma^{\pm 1}_n}{\alpha}_2
\]
where  ${\alpha}_1, {\alpha}_2 \in B_{m,n}.$

\noindent {\it (3)} \ An {\it algebraic $L$--move\/} is a parted $L$--move between algebraic mixed braids. From Fig.~16 one can easily derive analogous algebraic expressions for
an algebraic $L_o$--move and an algebraic $L_u$--move as in Eqs.~\ref{alglo} and \ref{alglu} for the classical situation, where now $\alpha_1, \alpha_2 \in B_{m,n}$ and $\alpha_1'$,
 $\alpha_2' \in B_{m,n+1}$ are obtained from $\alpha_1$,
 $\alpha_2$ by replacing each $\sigma_j$ by $\sigma_{j+1}$ for all
 $j=i,\ldots,n-1$.

\vspace{.05in}

{\it (4) \ Twisted loop conjugation} is defined to be the multiplication of an algebraic mixed braid by a loop $a_i$ or $a^{- 1}_i$ from the top and by its combing through the fixed braid $B$ from the bottom, and it has the algebraic expression:
\[
 \beta  \sim  {a^{\mp 1}_i} \beta {\rho^{\pm 1}_i}
\]
 for  $\beta, a_i, \rho_i \in B_{m,n}$,
where $\rho_i$ is the combing of the loop  $a_i$ through $B$. (Note that the
combing of the loop  ${a^{-1}_i}$ through $B$ is ${\rho^{-1}_i}$.)
\end{defn}

We stress that the  twisted loop conjugation is a move between algebraic mixed braids. For the first three moves view Figs.~16 and~17, where the braid $B$ is now the identity braid. We are now in a position to restate Theorem~\ref{partedcompl}  in terms of algebraic mixed braids.

\begin{thm}[Algebraic mixed braid equivalence for $S^3 \backslash \widehat B$] \label{algcompl}  Two  oriented
links in  $S^3 \backslash \widehat B$ are isotopic if and only if  any two
corresponding algebraic mixed braid representatives in  $B_{m,\infty}$  differ by a finite
sequence of algebraic braid relations and the following moves:

\noindent (1) \, algebraic stabilization moves,

\noindent (2) \, algebraic Markov conjugations,

\noindent (3) \, twisted loop conjugations.

Equivalently, by a finite sequence of algebraic braid relations and  the following moves:

\noindent (1$'$) \,  algebraic $L$--moves,

\noindent (2$'$) \, twisted loop conjugations.
\end{thm}

\subsection{Braid equivalence in c.c.o. $3$-manifolds}

Let $M=\chi (S^3, \widehat B)$. A link $L$ in $M$ may be seen as a link in $S^3 \backslash \widehat B$ with the extra freedom to slide across the 2--discs bounded in $M$ by the specified longitudes of the components of the surgery link  $\widehat B$.

\begin{defn} \rm  Let  $b$ be the oriented boundary of a ribbon  and let $L_1 \bigcup \widehat B$
and $L_2 \bigcup \widehat B$ be two oriented mixed links, so that $L_2 \bigcup \widehat B$ is
the band connected sum (over $b$) of a component, $c$, of $L_1$ and the specified (from the
framing) longitude of a surgery component of $\widehat B$.  This is not an isotopy between the mixed links in $S^3$ but it reflects isotopy between the links $L_1$ and  $L_2$ in $M$ and we shall call it a {\it band
move}.
\end{defn}

A band move can be, thus, performed in two steps: firstly, one of the small edges of $b$ is
glued to a part of $c$ so that the orientation of the band agrees with the orientation of  $c$, while the
other small edge of $b$, which we shall call {\it little band}, approaches a surgery component of $\widehat B$ in an arbitrary way. Secondly, the little
band is replaced by an arc running in parallel with the specified longitude of the surgery
component, such that the orientation of the arc agrees with the orientation of $b$ and
the resulting link is $L_2 \bigcup \widehat B$.

The first step of a band move is an isotopy move between the links $L_1$ and  $L_2$ in $S^3 \backslash \widehat B$. So, from now on whenever we say  {\it band move} we shall always be referring to the second step. Since $\widehat B$ is oriented in our setting, there are two types of
band moves according as in the second step the orientation of the arc replacing
the little band  agrees (type  $\alpha$)  or disagrees (type  $\beta$)  with the orientation of
the surgery component --and implicitly of its specified longitude. View Fig.~22, where $p$ is the integral framing of the surgery component. The discussion can be summarized by saying that two oriented links $L_1$, $L_2$ are isotopic in $M$ if and only if the mixed links $\widehat{B}\bigcup L_1$ and $\widehat{B}\bigcup L_2$ differ in $S^3$ by ambient isotopies that keep $\widehat{B}$ pointwise fixed and by the band moves (see \cite{LR1} for details).

\begin{figure}
     \begin{center}
     \includegraphics[width=14cm]{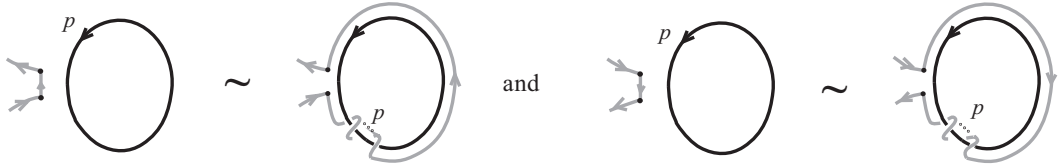}
     \caption{The two types of band moves }
     \label{}
\end{center}
\end{figure}

\smallbreak
Note now that neither of the two types of band moves can appear as a move between braids; so, in
order to extend our geometric mixed braid equivalence to $\chi (S^3, \widehat B)$, we modify appropriately the band move of type $\alpha$, by twisting the little band before performing the move resulting in a single braid band move that captures on the braid level both types of band moves. Namely:

\begin{defn} \rm A {\it braid band move}
is a move between mixed braids that upon closure is a band move. View Fig.~23, where only the move has been focused upon. The replacement of the little band links only with one of the strands of the same surgery component and runs in parallel to all remaining strands of that surgery component. A braid band move can be positive or negative depending on the type of crossing we choose for performing it.

\begin{figure}
     \begin{center}
     \includegraphics[width=8.5cm]{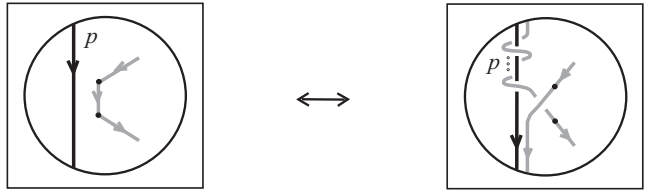}
     \caption{ The braid band move}
     \label{}
\end{center}
\end{figure}

\end{defn}

  We can now extend Theorem~\ref{lgeom} to the following result (cf. Theorem~5.10\cite{LR1}).

\begin{thm}[Geometric mixed braid equivalence for $\chi (S^3, \widehat B)$]\label{lgeomsurg} Let $L_1$, $L_2$ be two oriented links in $\chi (S^3, \widehat B)$ and let  $B_1 \bigcup B$,  $B_2 \bigcup B$  be two corresponding geometric mixed braids in $S^3$. Then  $L_1$  is isotopic to $L_2$  if and only if  $B_1 \bigcup B$  is equivalent to $B_2 \bigcup B$ in $S^3$ by a finite sequence of braid isotopies, the braid band moves and $L$--moves that do not affect $B$.
\end{thm}

We would like next to extend Theorem~\ref{partedcompl} to parted mixed braids related to c.c.o. $3$--manifolds. We first define:

\begin{defn}\rm  {\it A parted band move\/} is a braid band move
between parted mixed braids, such that: it takes place at the top of the braid (before any crossings of the surgery braid  are encountered), and the little band starts from the last strand of the moving
subbraid and moves over each strand of the parted mixed braid until it reaches
 the specified surgery strand from the right. Then, after the band move is performed, we apply to the resulting mixed braid the standard parting, bringing the new strands over to the last position of the moving subbraid. View  Fig.~24 for an example, where the moving subbraid
 has been simplified to the identity braid on two strands.
\end{defn}

\begin{figure}
\begin{center}
\includegraphics[width=7.5cm]{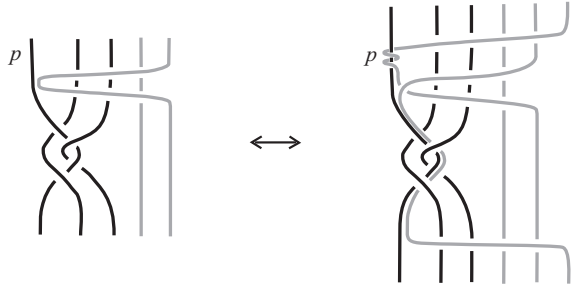}
\end{center}
\caption{  A parted band move }
\label{figure16}
\end{figure}

\begin{lem}[cf. Lemma~5\cite{LR2}] A braid band move may be always assumed, up to $L$--moves that do not affect $B$ and mixed braid isotopy, to take
place at the top part of a mixed braid and on the right of the specific surgery strand.
 Moreover, performing a braid band move on a parted mixed braid and then parting, the
result is equivalent, up to $L$--moves and loop conjugation, to performing a parted band
move.
\end{lem}

Lemma~5 enables us to extend Theorem~\ref{partedcompl} to parted mixed braids related to c.c.o. $3$--manifolds (cf. Theorem~3\cite{LR2}):

\begin{thm}[Parted mixed braid equivalence for $\chi (S^3, \widehat B)$]\label{partedcco} Two
oriented links in  $\chi (S^3, \widehat B)$ are isotopic if and only if any two
corresponding  parted mixed braids in $C_{m,\infty}$ differ by a finite sequence of parted mixed braid isotopies, parted
$L$--moves, loop  conjugations and parted band moves.

Equivalently, by a finite sequence of parted mixed braid isotopies, stabilization moves, Markov conjugations, loop conjugations and parted band moves.
\end{thm}

In order to move toward a completely algebraic statement we need
to understand  how a parted  band move is combed through the surgery braid $B$ and to
give algebraic expressions for parted band moves between algebraic mixed braids.

\begin{defn}\label{twistband} \rm  An  {\it algebraic  band move\/} is defined to be a  parted  band move
between elements of $B_{m,\infty}$. View Fig.~25 for an abstract example. Note that the isotopy of the
little band in the dotted box is treated as `invisible', that is, as identity in the braid
group. Setting:
\[
\lambda_{n-1} := \sigma_{n-1} \ldots \sigma_1  \mbox{ \ \ and \ \ }
t_{k,n} := \sigma_n \ldots \sigma_1 a_k {\sigma^{-1}_1} \ldots {\sigma^{-1}_n},
\]
an algebraic band move has the following algebraic expression:
\[
\beta_1 \beta_2 \ \sim \ \beta'_1 \, {t^{p_k}_{k,n}} \,
{\sigma^{\pm 1}_n} \, \beta'_2,
\]
\noindent where $\beta_1, \beta_2 \in B_{m,n}$ and $\beta'_1, \beta'_2  \in
B_{m,n+1}$  are the words $\beta_1, \beta_2$ respectively with the substitutions:
\[
\begin{array}{lcl}
  {a^{\pm 1}_k} & \longleftrightarrow & {[({\lambda^{-1}_{n-1}}
{\sigma^{2}_n} \lambda_{n-1}) \, a_k]}^{\pm 1}    \\

  {a^{\pm 1}_i} & \longleftrightarrow &   ({\lambda^{-1}_{n-1}} {\sigma^{2}_n}
 \lambda_{n-1}) \, {a^{\pm 1}_i} \, ({\lambda^{-1}_{n-1}} {\sigma^{2}_n}
\lambda^{-1}_{n-1}),   \mbox{ \ \ if \ } i < k   \\

 {a^{\pm 1}_i}  & \longleftrightarrow & {a^{\pm 1}_i},  \mbox{ \ \ if \ } i > k.    \\
 \end{array}
\]

Further, a {\it twisted algebraic band move\/} is defined to be a parted band move that has been combed through $B$, so it is the composition of an algebraic band move with the combing of the parallel strand, view  Fig.~26. Moreover, it has the following algebraic expression:
\[
\beta_1 \beta_2 \ \sim \ \beta'_1 \, {t^{p_k}_{k,n}} \,
{\sigma^{\pm 1}_n} \, \beta'_2 \, r_k
\]
\noindent where $r_k$ is the combing
 of the parted parallel strand to the $k$th surgery strand through the surgery braid.
 \end{defn}
We stress that a twisted algebraic band move is a move between algebraic mixed braids.

\begin{figure}
\begin{center}
\includegraphics[width=12cm]{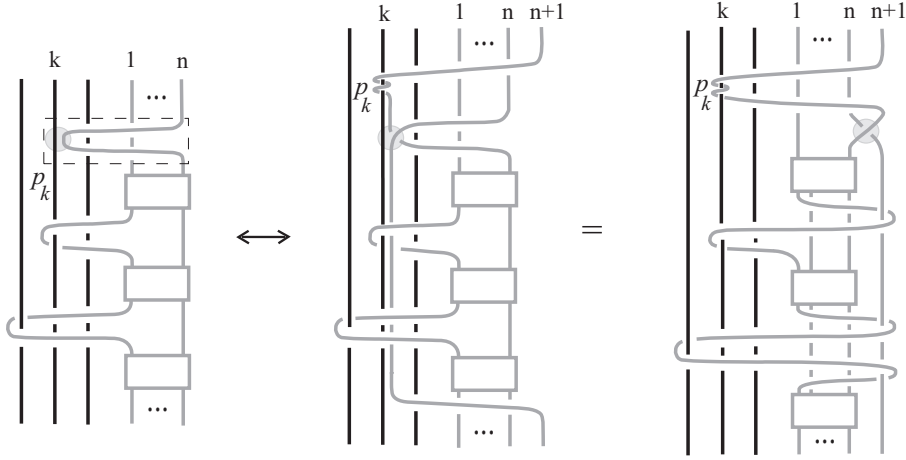}
\end{center}
\caption{  An algebraic band move and its algebraic expression }
\label{figure20}
\end{figure}

\begin{figure}
\begin{center}
\includegraphics[width=8cm]{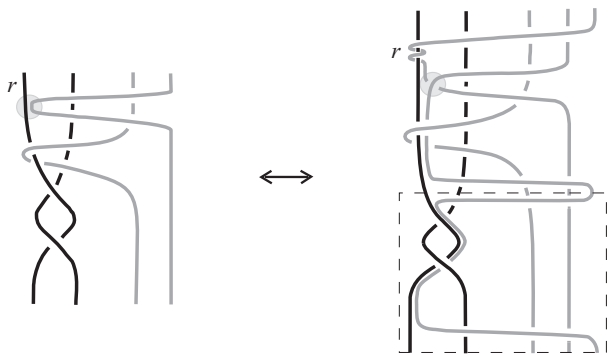}
\end{center}
\caption{  Parted band move =  algebraic band move + combing  }
\label{figure21}
\end{figure}

The above definitions lead immediately to the following (cf. Lemma~8\cite{LR2}):

\begin{lem}\label{combband}  Performing a parted band move on a parted mixed braid and then combing, the
result is the same as combing the mixed braid and then performing an algebraic band move.
\end{lem}

We are, finally, in the position to state the following result.

\begin{thm}[Algebraic mixed braid equivalence for $\chi (S^3, \widehat B)$]\label{algcco}  Two  oriented links
in  $\chi (S^3, \widehat B)$ are isotopic if and only if any two
corresponding algebraic mixed braid representatives in  $B_{m,\infty}$  differ by a finite
sequence of algebraic mixed braid relations and the following moves:

\noindent (1) \, algebraic stabilization moves (Definition~\ref{algmoves}),

\noindent (2) \, algebraic Markov conjugations (Definition~\ref{algmoves}),

\noindent (3) \, twisted loop conjugations (Definition~\ref{algmoves}),

\noindent (4) \,  twisted algebraic band moves (Definition~\ref{twistband}).

Equivalently, by a finite sequence of algebraic mixed braid relations and the following moves:

\noindent (1$'$) \,  algebraic $L$--moves (Definition~\ref{algmoves}),

\noindent (2$'$) \, twisted loop conjugations (Definition~\ref{algmoves}),

\noindent (3$'$) \ twisted algebraic band moves (Definition~\ref{twistband}).
\end{thm}

\begin{proof} By Theorem~\ref{algcompl} we only have to consider the case when a parted band
move takes place, and by Theorem~\ref{partedcco} we only have to consider the behaviour of a parted
band move with respect to combing, and this is Lemma~\ref{combband}.
\end{proof}

In \cite{LR2} we give the precise algebraic braid equivalence moves for knot complements and c.c.o. 3--manifolds defined by: the $m$--unlink, the Hopf link, a daisy chain, the right-handed trefoil.

\begin{rem}\rm The element ${t^{p_k}_{k,n}}$ in Definition~\ref{twistband} of an algebraic band
move is just a Markov conjugate of the loop ${a^{p_k}_k}$ and these are the appropriate words for defining inductive Markov traces on quotient algebras of group algebras of $B_{m,n}$.  Note also that the words in the parentheses of the substitutions of the loops get significantly simplified if we apply a quadratic relation on the
$\sigma_i$'s. Moreover, in Theorem~\ref{algcco} we obtain the most local description and the best possible control over the band moves of links in closed $3$--manifolds, and this is very useful for the study of skein modules of closed $3$--manifolds\cite{Pr}.
\end{rem}

\subsection{Braid equivalence in handlebodies}

A link inside a handlebody $H_m$ may not pass beyond the boundary of the handlebody from either end. Representing $H_m$ by the identity braid $I_m$ with infinitely extended strands, this means that a link inside $H_m$ may not cross the closing arc $k$ of Fig.~13. So, two oriented links $L_1, L_2$ in a handlebody $H_m$ are isotopic if and only if the mixed links $I_m \bigcup L_1$ and $I_m \bigcup L_2$ are isotopic in $S^3$ by an ambient
isotopy which keeps $I_m$ (with infinitely extended strands) pointwise fixed.

The fact that crossing the closing arc $k$ is {\it forbidden} in the category of handlebodies affects the definition of the closure of a geometric mixed braid in this category. Indeed, let $L$ be an oriented link in $H_m$ giving rise to the oriented mixed link $I_m \bigcup L$ in $S^3$, which gets braided to the geometric mixed braid $I_m \bigcup \beta$ (recall the braiding discussion for handlebodies in Subsection~2.1).  In order to ensure that the closure of $I_m \bigcup \beta$ will be isotopic to $I_m \bigcup L$, the {\it closure} is defined by joining with simple closing arcs the  endpoints of the corresponding strands of the moving part only, each closing arc passing {\it over} or {\it under} the mixed braid, according to the instruction `u' and `o' on the top endpoint of the geometric mixed braid, recall Fig.~9c. View Fig.~27 for an abstract illustration of the closure. To emphasize that different choices of closing labels will give rise, in general, to non-isotopic links in $H_m$ view the example in Fig.~28.

\begin{figure}
     \begin{center}
     \includegraphics[width=3.2cm]{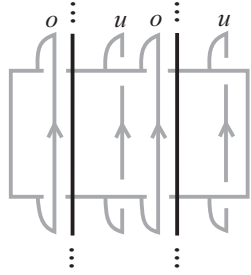}
     \caption{ The closure of a geometric mixed braid related to a handlebody }
     \label{}
\end{center}
\end{figure}

\begin{figure}
     \begin{center}
     \includegraphics[width=12cm]{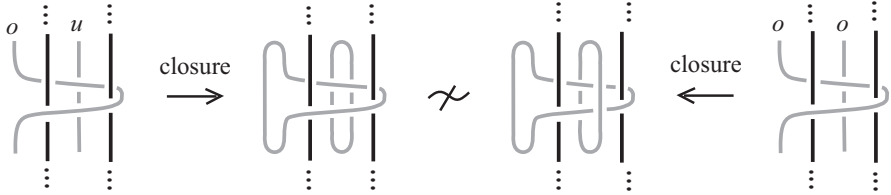}
     \caption{ Different closing labels yield non-isotopic links }
     \label{}
\end{center}
\end{figure}

The isotopy restriction for handlebodies, as supposed to knot complements, leads to the following definition of $L$--moves (cf. \cite{HL}):

\begin{defn}[Geometric $L$--moves in $H_m$] \rm Let $I_m \bigcup \beta$ be a braid in $H_m$. A {\it geometric $L$--move} consists in cutting a strand of the subbraid $\beta$ at a point, then pulling the two ends of the cut, the upper downward and the lower upward, both running entirely {\it over\/} the braid ({\it geometric $L_o$--move}) or entirely {\it under\/} it ({\it geometric $L_u$--move}), and keeping them aligned so as to obtain a  new pair of corresponding braid strands, which are labelled `o' and `u' according to the type of the $L$--move. View Fig.~29 for an example. Further, by small braid isotopies, an $L$--move is equivalent to adding an in-box crossing (positive or negative) facing the right-hand side (or, equally, the left-hand side) of the braid.
\end{defn}
Clearly, closing the new pair of strands created from an $L$--move results in a kink in the tangle diagram, so it corresponds to the isotopy move RI.

\begin{figure}
     \begin{center}
     \includegraphics[width=11.5cm]{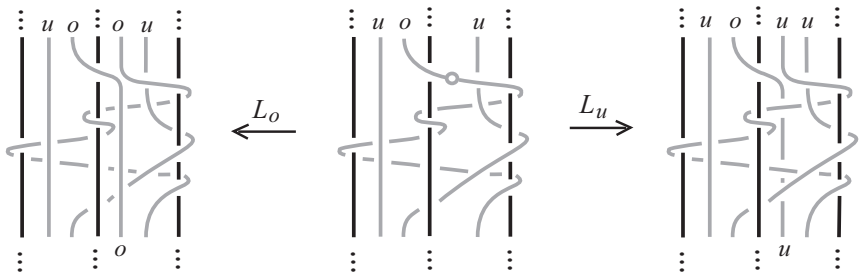}
     \caption{  The two types of $L$--moves in $H_m$ }
\end{center}
\end{figure}

We are now in a position to state the following (cf. Theorem~3 \cite{HL}).

\begin{thm}[Geometric braid equivalence for $H_m$]\label{geomhb}  Two oriented links in $H_m$ are isotopic if and only if any two corresponding geometric mixed braids differ by a finite sequence of geometric mixed braid isotopies and  geometric $L$--moves.
 \end{thm}

Moving toward an algebraic statement, we define the {\it parting} of a geometric mixed braid as described in Lemma~\ref{parting} and as illustrated in Fig.~14, only, here, the pulling of each pair of corresponding moving strands to the right is done according to its label `o' or `u'. Note that label `o' resp. `u' means to pull the pair of corresponding strands  over {\it all} resp. under {\it all} fixed strands to the right. So, if we would like to relate to the proof of Lemma~\ref{parting}, an array of labels here would contain the same symbol `o' or `u'.  View Fig.~30 for an example. A parted mixed braid related to a handlebody does not have labels `o' or `u' attached to its moving strands.

\begin{figure}
     \begin{center}
     \includegraphics[width=10cm]{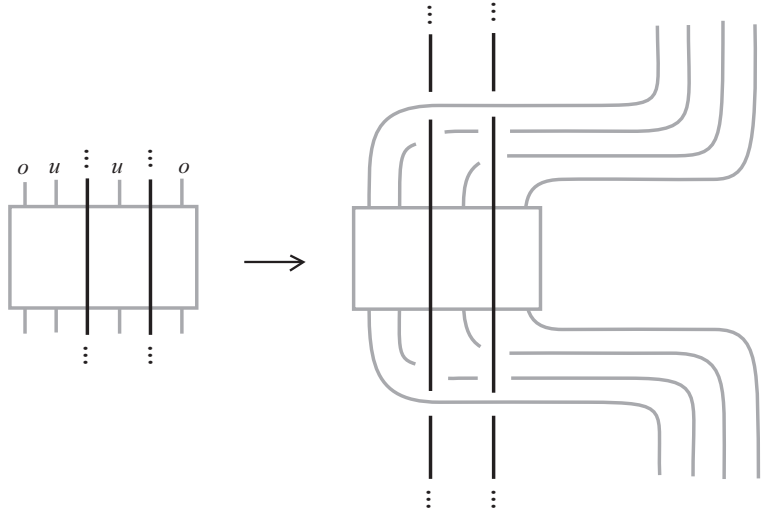}
     \caption{ The parting of a geometric mixed braid related to $H_m$ }
\end{center}
\end{figure}

Consequently, the {\it closure} of a parted mixed braid is realized by joining the pairs of corresponding endpoints of the moving part only with simple unlinked arcs, as in the standard closure of classical braids (view left-hand illustration of Fig.~31). This is not in contradiction with the fact that the parted mixed braids form a subset of the geometric mixed braids. Indeed, putting arbitrarily labels `o' and `u' to the endpoints of the moving strands of a parted mixed braid and applying the definition of closure of a geometric mixed braid will result in isotopic mixed links, as no surpassing of the fixed strands is involved (view right-hand illustration of Fig.~31).

\begin{figure}
     \begin{center}
     \includegraphics[width=8cm]{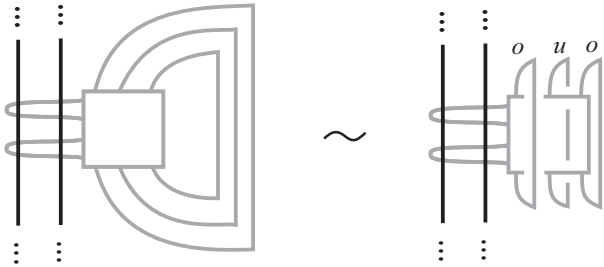}
     \caption{ The closure of a parted mixed braid }
\end{center}
\end{figure}

From the above, a parted mixed braid is the same as an algebraic mixed braid. So, we will be working with the braid groups $B_{m,n}$ and the algebraic moves of Definition~\ref{algmoves}. Note that, the fact that the parting cannot assume arbitrary labels (as it was the case for knot complements) forbids loop conjugation as a braid equivalence move for handlebodies.
Finally, the agreement of closing labels, $L$--moves and parting labels yields the following:

\begin{lem}  Parting a geometric mixed braid containing a geometric $L$--move gives rise to an algebraic mixed braid containing an algebraic $L$--move.
\end{lem}

The above together with Theorem~\ref{geomhb} lead to the following (cf. Theorems~4 and~5\cite{HL}):

\begin{thm}[Algebraic mixed braid equivalence for $H_m$]\label{alghb}  Two oriented links in
$H_m$ are isotopic if and only if any two corresponding algebraic mixed braids in  $B_{m,\infty}$ differ by a finite sequence of algebraic mixed braid relations and the following moves:

\noindent (1) \, algebraic stabilization moves (Definition~\ref{algmoves}),

\noindent (2) \, algebraic Markov conjugations (Definition~\ref{algmoves}).

Equivalently,  by a finite sequence of algebraic mixed braid relations and the following moves:

\noindent (1$'$) \,  algebraic $L$--moves (Definition~\ref{algmoves}).
\end{thm}

\section{Other diagrammatic settings}

 In this section we discuss braid equivalences for virtual knots, for flat virtuals, for welded knots and for singular knots. In all these settings the closure of the corresponding braids is realized as the ordinary closure of classical  braids.

\subsection{Virtual, flat virtual and welded braid equivalence }

 Virtual knot theory was introduced by Kauffman \cite{Kau2} and it is an extension of classical diagrammatic knot theory. In this extension one  adds  a {\em virtual crossing} that is neither an over-crossing nor an under-crossing. Virtual isotopy generalizes the ordinary Reidemeister
moves for classical links and is generated by planar isotopy and the moves RI, RII, RIII of Fig.~1 and the moves vRI, vRII, vRIII and the special detour move of Fig.~32, the last one being  the key move in the theory. Moves like F1 and F2 with two  real crossings and one virtual are forbidden in virtual knot theory. The virtual braid group $VB_{n}$ is generated by the classical crossings and the virtual crossings $v_i$, satisfying $v_i^2=1$, the braid relations in Eq.~\ref{brels} for the classical crossings and for the virtual crossings, and the mixed relations $v_i \sigma_{i+1} v_i = v_{i+1} \sigma_i v_{i+1}$ which are not symmetric. Cf. \cite{KL1} and references therein.

\begin{figure}
     \begin{center}
     \includegraphics[width=14cm]{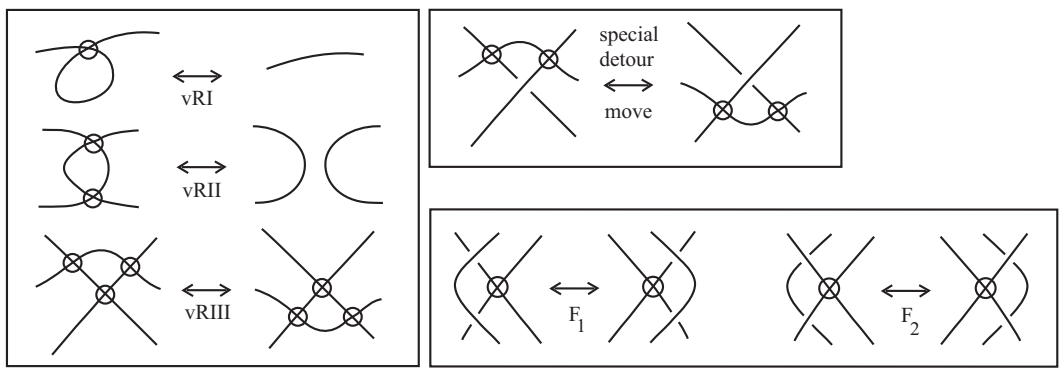}
     \caption{Virtual isotopy moves}
     \label{}
\end{center}
\end{figure}

A flat virtual knot is like a virtual knot but without the extra under/over structure at the real crossings. Instead we
have shadow crossings \cite{Tu}, the flat crossings. The study of flat virtual knots and links was initiated in
\cite{Kau2}. The isotopy moves and the two forbidden moves in virtual knot theory are completely analogous for the
flat virtual setting. In the flat virtual braid group on $n$ strands the generators  are the virtual crossings
$v_i$ and the flat crossings $c_i$, such that $c_i^2 = 1$. The  mixed relations  $v_{i} c_{i+1}
v_{i}  =  v_{i+1} c_i v_{i+1}$  are here also not symmetric. Cf. \cite{KL1} and references therein.

Welded  braids were introduced in \cite{FRR}. Welded knots satisfy the same isotopy relations as
the virtuals, but for welded knots  one of the two forbidden moves of Fig.~32 is allowed, the move F1 which
contains an over arc and one virtual crossing.  The welded braid group  on $n$ strands is a quotient of the
virtual braid group, so it can be presented with the same generators and relations, with the additional
relations: $v_{i} \sigma_{i+1} \sigma_{i}  =  \sigma_{i+1} \sigma_i v_{i+1}$  \ (F1).

Braid equivalence theorems for the above categories have been established in \cite{Ka} and \cite{KL2}. We will present here the $L$--move approach of \cite{KL2}.

\begin{defn}\rm
A  {\it basic $L_v$--move}  on a virtual braid,
consists in cutting an arc of the  braid open and pulling the upper cutpoint downward and
the lower  upward, so as to create a new pair of braid strands with corresponding endpoints, and such that both strands  cross entirely  {\it virtually}
 with the rest of the braid. See Fig.~33. In abstract illustrations this is indicated by placing virtual crossings on
the border of the braid box. See also Fig.~37 for a concrete example. Further, by a small braid isotopy that does not  change the relative positions of endpoints,
a basic $L_v$--move can be equivalently seen as introducing an in-box virtual crossing to a virtual braid
which faces either the {\it right} or the {\it left} side of the braid. If we want to
emphasize the existence of the virtual crossing, we will say {\it virtual $L_v$--move}, abbreviated to {\it
$vL_v$--move}. View  Fig.~33 for abstract illustrations.
\end{defn}

\begin{figure}
     \begin{center}
     \includegraphics[width=13.5cm]{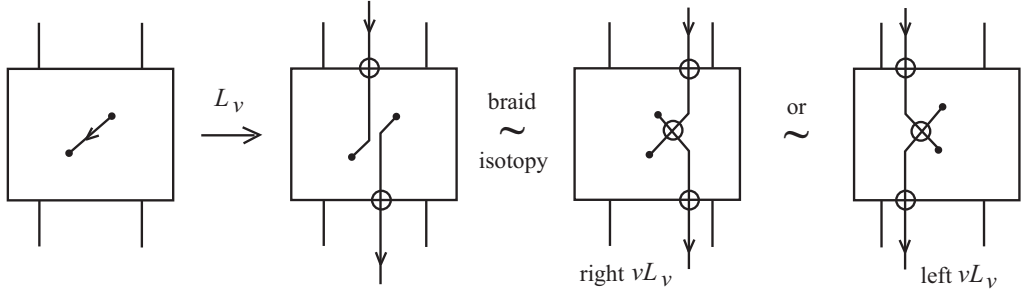}
     \caption{A basic $L_v$--move and the two $vL_v$--moves}
     \label{}
\end{center}
\end{figure}

In the closure of a basic $L_v$--move or a $vL_v$--move the detoured loop contracts to a virtual kink. But a kink
could  be also created by a real crossing, positive or negative. So we define:

\begin{defn}\rm
A  {\it real $L_v$--move}, abbreviated to  {\it  $+L_v$--move} or {\it  $-L_v$--move}, is a virtual
$L$--move that introduces  a real in-box crossing (positive or negative) on a virtual braid, and it can face either
the {\it right} or the {\it left} side of the braid. View Fig.~34 for abstract illustrations.
\end{defn}

\begin{figure}
     \begin{center}
     \includegraphics[width=12cm]{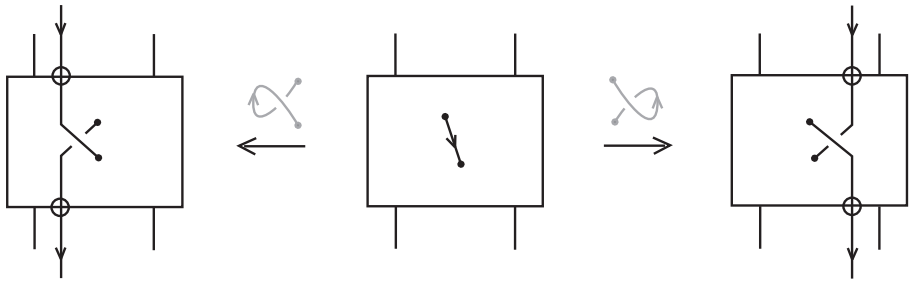}
     \caption{Left and right real $L_v$--moves}
     \label{}
\end{center}
\end{figure}

Further, if the crossing of the kink is virtual, then, in the presence of the forbidden moves, there
is another possibility for an $L_v$--move  on the braid level. Namely, we have:

\begin{defn}\rm
A  {\it threaded $L_v$--move} on a virtual braid is a virtual  $L$--move with a virtual
crossing in which, before pulling open the little up-arc of the kink, we perform a Reidemeister II
move with real crossings, using another arc of the braid, the {\it thread}. See Fig.~35. There are two possibilities:
an  {\it over-threaded $L_v$--move} and an {\it under-threaded $L_v$--move}, depending on whether we pull the kink
over or under the thread,  both with the variants {\it right} and {\it left}.
\end{defn}

\begin{figure}
     \begin{center}
     \includegraphics[width=11cm]{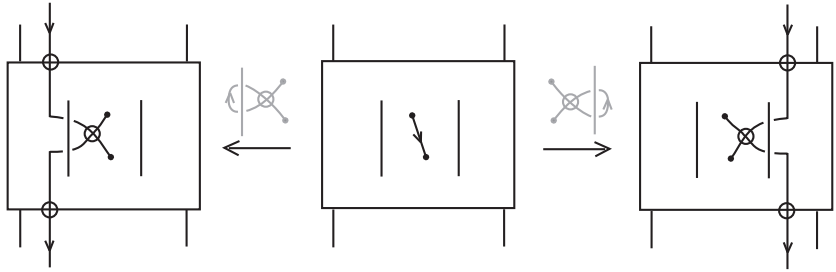}
     \caption{ Left and right under--threaded $L_v$--moves}
     \label{}
\end{center}
\end{figure}

A threaded $L_v$--move cannot be simplified in the braid.
 If the crossing of the kink were real, then, using a braided RIII move
with the thread, the move would reduce  to a real $L_v$--move. Similarly, if the
forbidden moves were allowed, a threaded $L_v$--move would reduce to a  $vL_v$--move.

The effect of a  virtual $L$--move, basic, real or threaded, is to
stretch (and cut open) an arc of the braid around the braid axis using the detour move, after twisting it and possibly
after threading it. Conversely, such a move between virtual braids gives rise to isotopic closures,
since the virtual $L$--moves shrink locally to kinks (grey diagrams in Figs.~34 and 35).
We may finally introduce the notion of a classical $L$--move.

\begin{defn}\rm \label{classicL}
  An {\it allowed classical $L_o$--move} resp.  {\it $L_u$--move } on a virtual braid
is like an $L$--move between classical braids (recall Definition~\ref{lmove}), such that the closures of the virtual braids before and after the move are isotopic. See
Fig.~36 for abstract illustrations.
An allowed  classical  $L$--move (over or under) may also introduce an in-box
crossing (positive, negative or virtual or it may even involve a thread, in the virtual setting).
\end{defn}

 In order that a classical $L$--move between virtual braids is allowed, in the sense that it gives rise to
isotopic virtual links upon closure, it is required that the virtual braid  has no virtual crossings on the
entire vertical zone either to the left or to the right of the new strands of the $L$--move, so that we can perform the closure of the new strands on that side. Otherwise, the presence of forbidden moves preclude such moves.
As it turns out, the allowed classical $L$--moves can be expressed in terms of
$L_v$--moves and real conjugation \cite{KL2}.

\begin{figure}
     \begin{center}
     \includegraphics[width=10cm]{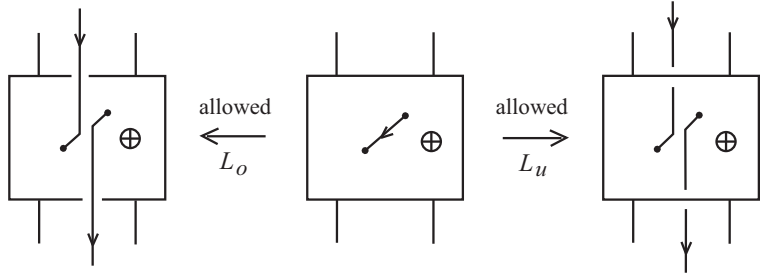}
     \caption{The allowed classical  $L$--moves}
     \label{}
\end{center}
\end{figure}

 In Fig.~37  we illustrate an example of various types of $L$--moves taking place at the same point
of a virtual  braid. We are now in a position to state the following results. For details consult \cite{KL2}.

\begin{figure}
     \begin{center}
     \includegraphics[width=12cm]{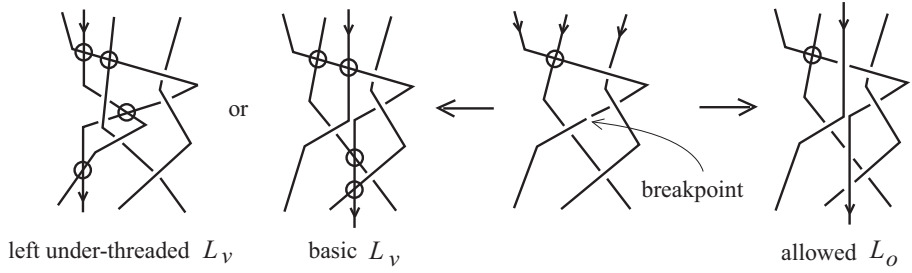}
     \caption{A concrete example of introducing $L$--moves in a virtual braid}
     \label{}
\end{center}
\end{figure}

\begin{thm}[Geometric equivalence for virtual braids]\label{geomv} Two oriented  virtual links are isotopic if and
only if any two corresponding virtual braids differ by virtual braid isotopy and a finite sequence of the
following moves:

\noindent (i) \ real conjugation,

\noindent (ii) \, right virtual $L_v$-moves,

\noindent (iii) \, right real $L_v$-moves,

\noindent (iv) \,  right and left under-threaded $L_v$-moves.
\end{thm}

In \cite{KL2} we stated the following conjecture.

 \begin{conj} \rm {\it Real conjugation is not a consequence of the $L_v$-moves.} In other words, it should be
possible to construct a virtual braid invariant  that will not distinguish  $L_v$-move equivalent virtual
braids, but will distinguish virtual braids that differ by real conjugation. As the simplest possible
puzzle, try to show that there is no sequence of $L_v$-moves connecting the  pair of equivalent
braids shown in Fig.~38.
\end{conj}

\begin{figure}
     \begin{center}
     \includegraphics[width=3cm]{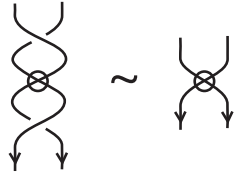}
     \caption{The simplest pair of real conjugates}
     \label{}
\end{center}
\end{figure}

\begin{thm}[Algebraic equivalence for virtual braids]\label{algv} Two oriented  virtual links are isotopic if and
only if any two corresponding virtual braids  in $\cup_nVB_n$ differ by a finite sequence of braid relations and the
following moves:

\noindent (i) \ virtual and real conjugation:    \  $ v_i \alpha v_i \sim \alpha \sim
{\sigma_i}^{-1}\alpha \sigma_i $,

\noindent (ii) \, right virtual and real stabilization: \  $\alpha v_n \sim \alpha
\sim \alpha \sigma_n^{\pm 1}$,

\noindent (iii) \, algebraic right under-threading (Fig.~39): \  $\alpha \sim \alpha \sigma_n^{-1} v_{n-1} \sigma_n^{+1} $,

\noindent (iv) \,  algebraic left under-threading (Fig.~39):  \  $\alpha \sim  \alpha v_n v_{n-1} \sigma_{n-1}^{+1} v_n
\sigma_{n-1}^{-1} v_{n-1} v_n $,

\noindent where $\alpha,  v_i, \sigma_i \in VB_n$ and  $v_n, \sigma_n \in VB_{n+1}$.
\end{thm}

\begin{figure}
     \begin{center}
     \includegraphics[width=10cm]{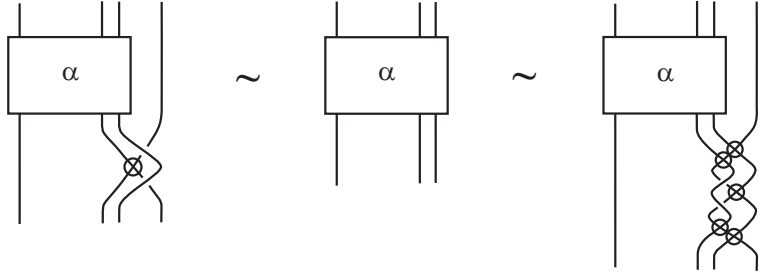}
     \caption{ Algebraic right and left under-threading}
     \label{}
\end{center}
\end{figure}

The statements for flat virtual braid equivalence are completely analogous. The only difference is that real crossings are substituted by flat crossings, cf.  \cite{KL2}. For the welded braid equivalence, move (iv) of Theorem~\ref{geomv} is not needed and moves (iii) and (iv) of Theorem~\ref{algv} are not needed, cf. \cite{KL2}.

We close by  pointing out that in \cite{KL2} it would be quite difficult to compare our local algebraic formulation of the braid equivalence for virtuals with that of Kamada \cite{Ka} without the fundamental $L$--move context.

\subsection{Singular braid equivalence}

Singular knots are related to Vassiliev's theory of knot invariants. In singular knot theory isotopy is generated by
planar isotopy, by the moves RI, RII, RIII of Fig.~1 and by the moves S1, S2, S3, S4 of Fig.~40. Moves SF1 and SF2 of Fig.~40 are forbidden in this theory, see \cite{Kau1}.

\begin{figure}
     \begin{center}
     \includegraphics[width=13cm]{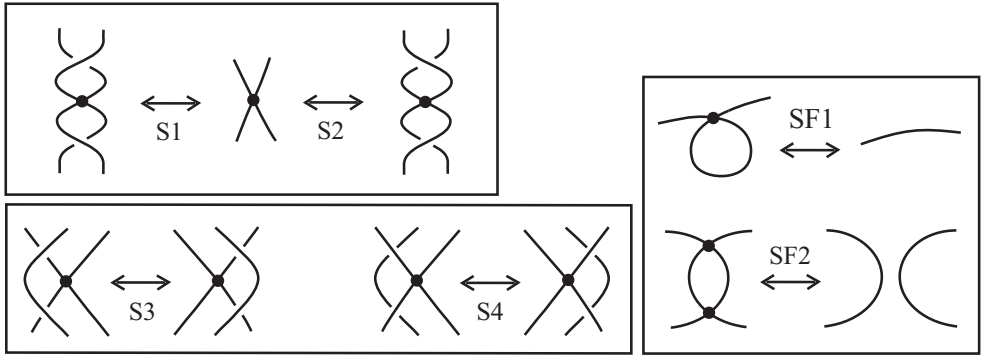}
     \caption{ Diagrammatic moves for singular knots }
     \label{}
\end{center}
\end{figure}

The singular crossings $\tau_i$ generate, together with the real crossings $\sigma_i$ and their inverses, the
singular braid monoid $SB_n$ (introduced in different contexts by Baez\cite{Ba}, Birman\cite{Bi2} and Smolin\cite{Sm})  and they satisfy the braid relations Eq.~\ref{brels} and the relations below.
\begin{equation}\label{relsbn}
\begin{array}{rclcll}
\sigma_i\sigma_i^{-1}    & = & \sigma_i^{-1} \sigma_i   & = & 1 &  \mbox{for all}\,\, i \\
\sigma_{i+1} \sigma_i \tau_{i+1} & = & \tau_i \sigma_{i+1} \sigma_i &  &  & \mbox{for}\, \vert i-j\vert =1 \\
 \sigma_i \tau_j & = & \tau_j\sigma_i &  &  & \mbox{for}\, \vert i-j\vert  >1 \\
 \sigma_i \tau_i & = & \tau_i\sigma_i & & & \mbox{for all}\,\, i \\
 \tau_i \tau_j & = & \tau_j\tau_i &  &  & \mbox{for}\, \vert i-j\vert  >1 \\
\end{array}
\end{equation}

The algebraic braid equivalence for singular braids was established by Gemein \cite{Ge}:

\begin{thm}[Algebraic equivalence for singular braids \cite{Ge}]\label{algs}
 Two singular braids in $\cup_n SB_n$ have isotopic closures if and only if they differ by a finite sequence of braid relations in $\cup_n SB_n$ and the following moves:

\noindent (i) \, singular commuting: \ $\tau_i\omega \sim \omega \tau_i, \quad \omega,\tau_i \in SB_n$

\noindent (ii) \ real conjugation: \ $\sigma_i \omega \sim \omega \sigma_i, \quad \omega,\sigma_i \in SB_n$

\noindent (iii) real stabilization: \ $\omega \sim \omega \sigma_n^{\pm 1}, \quad \omega \in SB_n$.
\end{thm}

Using the above result, in \cite{La3} the geometric braid equivalence for singular braids was easily established:

\begin{thm}[Geometric equivalence for singular braids] Two oriented  singular links are isotopic if and
only if any two corresponding singular braids differ by a finite
sequence of braid relations in $\cup_n SB_n$ and the following moves:

\noindent (i) singular commuting: \ $\tau_i\omega \sim \omega \tau_i, \quad \omega,\tau_i \in SB_n$

\noindent (ii) classical $L$--moves (Definition~\ref{lmove}).
\end{thm}



\begin{thebibliography}{99}

\bibitem[Al]{Al} J. W. Alexander,  A lemma on systems of knotted curves,
{\it Proc. Nat. Acad. Sci. U.S.A.} {\bf 9} (1923) 93--95.

\bibitem[Ar1]{Ar1} E. Artin, Theorie der Z\"{o}pfe, {\it Abh. Math. Sem.
              Hamburg Univ.} {\bf 4} (1925) 47--72.

\bibitem[Ar2]{Ar2} E. Artin, Theory of braids, {\it Ann. of Math.} (2) {\bf 48} (1947)
 101--126.

\bibitem[Be]{Ben} D. Bennequin, Entrlacements et \'{e}quations de
Pfaffe, {\it Asterisque} {\bf 107-108} (1983) 87--161.

\bibitem[Ba]{Ba} J. Baez,  Link invariants and perturbation theory, {\it Lett. Math. Phys.} {\bf 2} (1992) 43--51.

\bibitem[Bi1]{Bi1}  J. S. Birman, {\it Braids, links and mapping class groups},  Annals of Mathematics
             Studies, Vol. 82 (Princeton University Press, Princeton, 1974).

\bibitem[Bi2]{Bi2} J. Birman,  New points of view in knot theory, {\it Bull. Amer. Math. Soc. (N.S.)} {\bf 28}(2) (1993)
              253--287.

\bibitem[BM]{BM} J. S. Birman, W. W. Menasco, On Markov's Theorem, {\it J. Knot Theory and
             Ramifications} {\bf 11} no.3 (2002), 295--310.

\bibitem[Br]{Br} H. Brunn,  \"{U}ber verknotete Curven,  {\it Verh. des intern. Math.
             Congr.} {\bf 1}, 256--259 (1897).

\bibitem[Ch]{Ch} W.-L. Chow, On the algebraical braid group, {\it Ann. of Math.} (2) {\bf 49} (1948) 654--658.

\bibitem[FRR]{FRR} R. Fenn, R. Rimanyi, C.P. Rourke, The braid permutation group, {\em Topology} {\bf 36} (1997) 123--135.

\bibitem[Ge]{Ge}  B. Gemein,  Singular braids and Markov's theorem, {\em J. Knot Theory
              Ramifications} {\bf 6}(4) (1997) 441--454.

\bibitem[HL]{HL} R.~H\"{a}ring-Oldenburg, S.~Lambropoulou, Knot theory in handlebodies,
           {\it J. Knot Theory Ramifications} {\bf 11}(6) (2002) 921--943.

\bibitem[Jo]{Jo} V.~F.~R.~Jones, Hecke algebra representations of braid groups and link
               polynomials, {\it Annals of Math.} {\bf 126} (1987) 335--388.

\bibitem[Ka]{Ka} S. Kamada, Braid representation of virtual knots and welded knots,
            {\it  Osaka J. Math.} {\bf 44} (2007), no. 2, 441--458. See also {\em arXiv:math.GT/0008092}.

\bibitem[KT]{KT}  C. Kassel and V. Turaev, {\it Braid groups},  Graduate Texts in Mathematics, Vol. 247 (Springer, 2008).

\bibitem[Kau1]{Kau1}  L. H. Kauffman,  Invariants of graphs in three-space, {\em Trans. Amer. Math. Soc.} {\bf 311}(2) (1989)
              697--710.

\bibitem[Kau2]{Kau2} L. H. Kauffman, Virtual Knot Theory, {\em European J. Comb.} {\bf 20} (1999) 663--690.

\bibitem[KL1]{KL1}  L. H. Kauffman and S. Lambropoulou,  Virtual Braids, {\em Fundamenta Mathematicae} {\bf 184} (2004)
              159--186.

\bibitem[KL2]{KL2}  L. H. Kauffman and S. Lambropoulou,  Virtual Braids and the $L$-Move, {\em J. Knot Theory
               Ramifications} {\bf 15}(6) (2006) 773--811.

\bibitem[La1]{La1} {S. Lambropoulou}, ``A study of braids in 3--manifolds'',
             Ph.D. thesis, Warwick, 1993.

\bibitem[La2]{La2} S. Lambropoulou, Braid structures in handlebodies, knot complements and
              $3$-manifolds, in {\it Proceedings of Knots in Hellas '98}, Series on Knots and Everything Vol. 24
              (World Scientific, 2000), pp. 274--289.

\bibitem[La3]{La3} S. Lambropoulou, L-moves and Markov theorems, {\it J. Knot Theory Ramifications} {\bf 16} no. 10 (2007), 1--10.

\bibitem[La4]{La4}  S. Lambropoulou, ``Braid structures and braid equivalence in different manifolds and in different
              settings'', World Scientific, Series on Knots and Everything, in preparation.

\bibitem[LR1]{LR1} S. Lambropoulou and C. P. Rourke, Markov's theorem in three--manifolds, {\em Topology and its
               Applications} {\bf 78} (1997) 95--122.

\bibitem[LR2]{LR2} S. Lambropoulou and C. P. Rourke, Algebraic Markov equivalence for links in $3$-manifolds,
              {\em Compositio Math.} {\bf 142} (2006) 1039--1062.

\bibitem[Ma]{Ma} {A.A. Markov}, \"{U}ber die freie \"{A}quivalenz der
              geschlossenen Z\"{o}pfe, {\it Recueil Math\'{e}matique Moscou} {\bf 1}(43) (1936) 73--78.

\bibitem[Mo]{Mo} {H.R. Morton}, Threading knot diagrams,
                 {\it Math. Proc. Cambridge Philos. Soc.} {\bf 99} (1986) 247--260.

\bibitem[Pr]{Pr}{J.H. Przytycki},  Skein Modules, eprint: math.GT/0602264;
Chapter IX of the book ``KNOTS: From combinatorics of knot diagrams to the
combinatorial topology based on knots'', Cambridge University Press,to appear 2014, pp. 600.

\bibitem[Rd]{Rd1} {K. Reidemeister}, Elementare Begr\"{u}ndung der Knotentheorie, {\it Abh. Math. Sem.
              Hamburg Univ.} {\bf 5} (1927) 24--32.

\bibitem[Sk]{Sk} {R. Skora},  Closed braids in 3--manifolds,
             {\it Math. Zeitschrift} {\bf 211} (1992) 173--187.

\bibitem[Sm]{Sm} L. Smolin, {\em Knot theory, loop space and the diffeomorphism group}, New perspectives in canonical gravity,  245--266, Monogr. Textbooks Phys. Sci. Lecture Notes, {\bf 5}, Bibliopolis, Naples, 1988.

\bibitem[Su]{Su} {P.A. Sundheim},  Reidemeister's theorem for 3--manifolds,
              {\it Math. Proc. Camb. Phil. Soc.} {\bf 110} (1991)  281--292.

\bibitem[Tu]{Tu} {V. Turaev}, Virtual strings, {\em Ann. Inst. Fourier (Grenoble)}  {\bf 54}(7) (2004) 2455--2525.

\bibitem[Tr]{Tr} {P. Traczyk}, A new proof of Markov's braid theorem, preprint  (1992), Banach Center
Publications, Vol.42, Institute of Mathematics Polish Academy of Sciences, Warszawa (1998).

\bibitem[Vo]{Vo} {P. Vogel}, Representation of links by braids: A new algorithm, {\it Comment. Math.
             Helvetici} {\bf 65} (1990) 104--113.

\bibitem[We]{We}  {N. Weinberg},  Sur l' equivalence libre des tresses
             ferm\'{e}e, {\it Comptes Rendus (Doklady) de l' Acad\'{e}mie des Sciences de l' URSS}
              {\bf 23}(3) (1939)  215--216.

\bibitem[Ya]{Ya} {S. Yamada}, The minimal number of Seifert circles equals the braid index of a link,
             {\it Invent. Math.} {\bf 89} (1987) 347--356.


\end{thebibliography}
\end{document}